\documentclass{elsart}
\usepackage{amssymb}
\journal{Finite Fields and Their Applications}
\newtheorem{theorem}{Theorem}[section]
\newtheorem{corollary}[theorem]{Corollary}
\newtheorem{lemma}[theorem]{Lemma}
\newtheorem{proposition}[theorem]{Proposition}

\newcommand{\N}{\mathbb N}
\newcommand{\Z}{\mathbb Z}
\newcommand{\F}{\mathbb F}

\newcommand{\fq}{\F_{\hskip-0.7mm q}}
\newcommand{\cfq}{\overline{\F}_{\hskip-0.7mm q}}

\def\ifm#1#2{\relax \ifmmode#1\else#2\fi}

\newcommand{\ds}{\displaystyle}

\newcommand{\klk}    {\ifm {,\ldots,} {$,\ldots,$}}

\newcommand{\plp}    {\ifm {+\cdots+} {$+\ldots+$}}

\newcommand{\A}{\mathbb A}

\newcommand{\om}[2]   {{#1}_1 \klk {#1}_{#2}}

\newcommand{\xo}[1]  {\ifm {\om X {#1}} {$\om X {#1}$}}
\newcommand{\xon}    {\ifm {\om X n} {$\om X n$}}
\newcommand{\yo}[1]  {\ifm {\om Y {#1}} {$\om Y {#1}$}}

\newenvironment{proof}
 {\medskip\noindent {\sc Proof}. \ }
 {\hfill\vbox{\hrule height 5pt width 5pt } \bigskip}

\begin{document}
\begin{frontmatter}
\title{Improved  explicit estimates on the number of solutions of
equations over a finite field\thanksref{grants}}
\author[Cafure,ungs]{A. Cafure}\ead{acafure@dm.uba.ar},
\author[ungs,conicet]{G. Matera\corauthref{cor}}\ead{gmatera@ungs.edu.ar}
\ead[url]{www.medicis.polytechnique.fr/\~{}matera}

\address[Cafure]{Departamento de Matem\'atica, Facultad de Ciencias
Exactas y Naturales, Universidad de Buenos Aires, Ciudad
Universitaria, Pabell\'on I \newline (1428) Buenos Aires,
Argentina.}

\address[ungs]{Instituto de Desarrollo Humano,
Universidad Nacional de General Sarmiento, J.M. Guti\'errez 1150
(1613) Los Polvorines, Buenos Aires, Argentina.}
\address[conicet]{Member of the CONICET, Argentina.}

\corauth[cor]{Corresponding author.}

\thanks[grants]{Research was partially supported by the
following Argentinian and German grants: UBACyT X198, PIP CONICET
2461, BMBF--SETCIP AL/PA/01--EIII/02, UNGS 30/3005.  Some of the
results presented here were first announced at the {\em Workshop
Argentino de Inform\'atica Te\'orica}, WAIT'02, held in September
2002 (see \cite{CaMa02}).}

\begin{abstract}
We show explicit estimates on the number of $q$--rational points
of an $\fq$--definable affine absolutely irreducible variety of
$\cfq^n$. Our estimates for a hypersurface significantly improve
previous estimates of W. Schmidt and M.-D. Huang \& Y.-C. Wong,
while in the case of a variety our estimates improve those of S.
Ghorpade \& G. Lachaud in several important cases. Our proofs rely
on elementary methods of effective elimination theory and suitable
effective versions of the first Bertini theorem.
\end{abstract}
\end{frontmatter}
%
%
\section{Introduction.}
Let $p$ be a prime number, let $q:=p^k$, let $\fq$ denote the
finite field of $q$ elements and let $\cfq$ denote the algebraic
closure of the field $\fq$. Let be given a finite set of
polynomials $F_1,\dots,F_m\in\fq[X_1,\dots,X_n]$ and let $V$
denote the affine subvariety of $\cfq^n$ defined by
$F_1,\dots,F_m$. Counting or estimating the number of
$q$--rational points $x\in\mathbb{F}_q^n$ of $V$ is an important
subject of mathematics and computer science, with many
applications.\smallskip

In a fundamental work \cite{Weil48}, A. Weil showed that for any
$\fq$--definable absolutely irreducible plane curve $\mathcal{C}$
of degree $\delta$ and genus $g$, the following estimate holds:
$$|\#(\mathcal{C}\cap \fq^2)-q|\leq 2gq^{1/2}+\delta+1.$$
Taking into account the well--known inequality $2g \leq
(\delta-1)(\delta-2)$, we have the estimate
\begin{equation}\label{estimateWeil}
|\#(\mathcal{C}\cap \fq^2)-q|\leq (\delta-1)(\delta-2)q^{1/2} +
\delta+1,
\end{equation}
which is optimal in the general case. The proof of this result was
based on  sophisticated techniques of algebraic geometry.

Weil's estimate (\ref{estimateWeil}) was generalized to higher
dimensional varieties by S. Lang and Weil \cite{LaWe54}. Their
result may be rephrased as follows: for any $\fq$--definable
absolutely irreducible subvariety $V$ of $\cfq^n$ of dimension
$r>0$ and degree $\delta>0$, we have the estimate:
\begin{equation}\label{estimate:LangWeil}
|\#(V\cap \fq^n)-q^{r-1}|\leq (\delta-1)(\delta-2)q^{r-1/2}+
Cq^{r-1},\end{equation}
where $C$ is a universal constant, depending only on $n$, $r$ and
$\delta$, which was not explicitly estimated.

S. Ghorpade and G. Lachaud (\cite{GhLa02}, \cite{GhLa02a}) found
an explicit estimate on the constant $C$ of
(\ref{estimate:LangWeil}). More precisely, in \cite[Remark
11.3]{GhLa02a} (see also \cite[Theorem 4.1]{GhLa02}) the following
estimate on the number $N$ of $q$--rational points of an
$\fq$--definable absolutely irreducible subvariety $V$ of $\cfq^n$
of dimension $r>0$ and degree $\delta>0$ is shown:
\begin{equation}\label{estimate:Ghorpade-Lachaud}|N-q^r| \leq
(\delta-1)(\delta-2)q^{r-1/2}+6\cdot2^s(sd+3)^{n+1}q^{r-1},
\end{equation}
where $s$ is the number of equations defining the variety $V$ and
$d$ is an upper bound of the degrees of these equations. Observe
that in the case of a hypersurface, estimate
(\ref{estimate:Ghorpade-Lachaud}) gives
\begin{equation}\label{Ghorpadehypersurface}|N-q^{n-1}|\le
(\delta-1)(\delta-2)q^{n-1/2}+12(\delta+3)^{n+1}q^{n-2}.
\end{equation}
The proof of this result is based on a sophisticated method
relying on a generalization of the Weak Lefschetz Theorem to
singular varieties and estimates of the Betti numbers of suitable
spaces of \'etale $\ell$--adic cohomology.

On the other hand, the first general estimate obtained by
elementary means was given by W. Schmidt in \cite{Schmidt73} (see
also \cite{Schmidt76}, \cite{Bombieri73}, \cite{LiNi83}).
Generalizing a method of S. Stepanov \cite{Stepanov71}, Schmidt
obtained the estimate
$$ |\#(\mathcal{C}\cap \fq^2)-q|\leq
\sqrt{2}\delta^{5/2}q^{1/2},$$
where $\mathcal{C}\subset\cfq^2$ is an absolutely irreducible
$\fq$-- definable plane curve of degree $\delta>0$ and the {\em
regularity condition} $q > 250\delta^5$ holds.

Later on, using an adaptation of Stepanov's method to the
hypersurface case Schmidt \cite{Schmidt74} showed the following
non--trivial lower bound for any absolutely irreducible
$\fq\!$--definable hypersurface $H$ of degree $\delta>0$:
\begin{equation}\label{estimate:Schmidt74}
\#(H\cap\fq^n)>q^{n-1}-(\delta-1)(\delta-2)q^{n-3/2}-(5\delta^2+
\delta+1)q^{n-2},
\end{equation}
provided that the regularity condition
$q>cn^3\delta^5\log^3\delta$ holds for a suitable constant $c>0$.
Let us remark that, up to now, this was the best explicit lower
bound known for an arbitrary absolutely irreducible
$\fq$--hypersurface. He also obtained in \cite{Schmidt76} the
following explicit estimate:
$$|\#(H\cap\fq^n)-q^{n-1}|\leq (\delta-1)(\delta-2)q^{n-3/2} +
6\delta^2\theta^{2^{\theta}}q^{n-2},$$
with $\theta:={(\delta+1)\delta/2}$.

Finally, combining (\ref{estimate:Schmidt74}) with Schmidt's
\cite{Schmidt76} method and Kaltofen's effective version of the
first Bertini Theorem \cite{Kaltofen95a}, M.-D. Huang and Y.-C.
Wong \cite{HuWo98} obtained the following estimate for an
$\fq$--definable absolutely irreducible hypersurface $H$ of
$\cfq^n$ of degree $\delta$: 
\begin{equation}\label{estimate:HuangWong}
\!\!\!\!\!|\#(H\cap\fq^n)-q^{n-1}|\le
(\delta-1)(\delta-2)q^{n-3/2}
+(\delta^2+2\delta^5)q^{n-2}+2\delta^7q^{n-5/2},
\end{equation}
provided that the regularity condition
$q>cn^3\delta^5\log^3\delta$ holds.

From the point of view of practical applications it is important
to improve as much as possible the regularity condition underlying
(\ref{estimate:Schmidt74}) and (\ref{estimate:HuangWong}).
Furthermore, estimate (\ref{estimate:HuangWong}) may grow large in
concrete cases due to the powers of $\delta$ arising in the
right--hand side of (\ref{estimate:HuangWong}). This is also the
case of (\ref{estimate:Ghorpade-Lachaud}) and
(\ref{Ghorpadehypersurface}), whose right--hand sides include
terms which depend exponentially on $n$ and the number $s$ of
equations.

In this article, combining techniques of \cite{Schmidt74},
\cite{Schmidt76} and \cite{Kaltofen95a} we obtain improved
explicit estimates on the number of $q$--rational points of an
$\fq$--definable affine absolutely irreducible variety $V$ of
$\cfq^n$. Our estimates in the case of a hypersurface
significantly improve the regularity of (\ref{estimate:Schmidt74})
and extend it, providing a corresponding upper bound. Further, we
improve both the regularity  and the right--hand side of
(\ref{estimate:HuangWong}) and  exponentially improve
(\ref{Ghorpadehypersurface}). Finally, in the case of an
absolutely irreducible variety, the worst case of our estimates
improve (\ref{estimate:Ghorpade-Lachaud}) in several important
cases, such as those of low codimension (for example $2r \geq
n-1$) and those of low degree (for example $d \leq 2(n-r)$).

Our methods rely on elementary arguments of effective elimination
theory (see Sections \ref{section: Notions} and
\ref{section:reduction}). In particular, we obtain elementary
upper bounds on the number of $q$--rational points of certain
$\fq$--definable affine varieties which improve \cite{Schmidt74},
\cite{Schmidt76} and \cite{ChRo96} (see Section \ref{section:
Notions}).

Our estimate for a hypersurface combines ideas of \cite{Schmidt74}
with an effective version of the first Bertini Theorem due to E.
Kaltofen \cite{Kaltofen95a}. Kaltofen's result is based on the
analysis of an algorithm that decides whether a given bivariate
polynomial with coefficients in a field is absolutely irreducible.
In Section \ref{section:algoritmoygeneric}, we adapt Kaltofen's
algorithm in order to determine the existence of irreducible
factors of a given degree of the restriction of a multivariate
absolutely irreducible polynomial to a plane. This allows us to
obtain suitable upper bounds on the genericity condition
underlying the choice of a restriction having no irreducible
factors of a given degree (Theorem \ref{nuestrokaltofen}). In
Section \ref{section:number of planes} we combine this result with
a combinatorial approach inspired in \cite{Schmidt74} in order to
estimate the number of restrictions of a given absolutely
irreducible polynomial $f\in \fq[\xon]$ to affine planes having a
fixed number of absolutely irreducible factors over $\fq$.

In Section \ref{section:hypersurface}, applying the estimates of
the preceding Section and adapting the methods of
\cite{Schmidt76}, we obtain the following estimate for an
absolutely irreducible $\fq$--hypersurface $H\subset\fq^n$ of
degree $\delta$ (see Theorem \ref{theorem:est_pol_abs_irr}), which
holds without any regularity condition:
%
%
$$ |\#(H\cap \fq^n)-q^{n-1}|\le(\delta-1)(\delta-2)q^{n-1/2}
+5\delta^{\frac{13}{3}}q^{n-2}. $$
%
%
Furthermore, using the lower bound underlying the previous
estimate we obtain the following estimate (see Theorem
\ref{theorem:lower bound}): for $q
>15\delta^{\frac{13}{3}}$ we have
%
$$ |\#(H\cap \fq^n)-q^{n-1}|\leq
(\delta-1)(\delta-2)q^{n-1/2}+ (5\delta^2+\delta+1)q^{n-2}.$$

Finally, in Section \ref{section:reduction} we combine these
estimates with elementary methods of effective elimination theory
(see Propositions \ref{proposition:chow} and
\ref{proposition:morph_bir}) in order to obtain estimates for
affine $\fq$--definable affine varieties (see Theorems
\ref{theorem:estimate_hyper_gral}, \ref{theorem:est_final_abs_irr}
and \ref{theorem:est_final_variety}). As an illustration of the
results we obtain, we have the following 
estimate for an $\fq$--definable absolutely irreducible variety
$V\subset\cfq^n$ of dimension $r>0$ and degree $\delta$: for
$q>2(r+1)\delta^2$, there holds $$|\#(V\cap \fq^n)-q^r| \leq
(\delta-1)(\delta-2)q^{r-1/2}+ 5\delta^{\frac{13}{3}} q^{r-1}.$$
%
\section{Notions and notations.}\label{section: Notions}
We use standard notions and notations of commutative algebra and
algebraic geometry as can be found in e.g. \cite{Kunz85},
\cite{Shafarevich84}, \cite{Matsumura80}.

For a given $m\in\N$, we denote by $\A^m=\A^m(\cfq)$ the
$m$--dimensional affine space $\cfq^m$ endowed with the Zariski
topology.

Let $\xon$ be indeterminates over $\fq$ and let $\fq[\xon]$ be the
ring of $n$--variate polynomials in the indeterminates $\xon$ and
coefficients in $\fq$. Let $V$ be an $\fq$--definable affine
subvariety $V$ of $\A^n$ (an $\fq$--variety for short). We shall
denote by $I(V)\subset\fq[\xon]$ its defining ideal and by
$\fq[V]$ its coordinate ring, namely, the quotient ring
$\fq[V]:=\fq[\xon]/I(V)$.

If $V$ is irreducible as an $\fq$--variety ($\fq$--irreducible for
short), we define its {\em degree} as the maximum number of points
lying in the intersection of $V$ with an affine linear subspace
$L$ of $\A^n$ of codimension $\dim(V)$ for which $\#(V\cap
L)<\infty$ holds. More generally, if $V=C_1\cup\cdots\cup C_h$ is
the decomposition of $V$ into irreducible $\fq$--components, we
define the degree of $V$ as $\deg(V):=\sum_{i=1}^h\deg(C_i)$ (cf.
\cite{Heintz83}). In the sequel we shall make use of the following
{\em B\'ezout inequality} (see \cite{Heintz83}, \cite{Fulton84}):
if $V$ and $W$ are $\fq$--subvarieties of $\A^n$, then
\begin{equation}\label{equation:Bezout}\deg (V\cap W)\le \deg V
\deg W.
\end{equation}
An $\fq$--variety $V \subset\A^n$ is {\em absolutely irreducible}
if it is irreducible as $\cfq$--variety.
%
%
\subsection{Some elementary upper bounds.}
In this section we exhibit upper bounds on the number of
$q$--rational points of certain $\fq$--varieties using elementary
arguments of effective elimination theory and the B\'ezout
inequality (\ref{equation:Bezout}). The purpose of this section is
to illustrate how these arguments significantly simplify the
previous combinatorial proofs (cf. \cite{LiNi83},
\cite{Schmidt74}, \cite{Schmidt76}), yielding also better
estimates than the usual ones in some cases. We start with the
following well--known result:
\begin{lemma}\label{lemma:hesch2} Let $V\subset\A^n$ be
an $\fq$--variety of dimension $r\ge 0$ and degree $\delta>0$.
Then the inequality $\#(V\cap\fq^n)\le \delta q^{r}$ holds.
\end{lemma}
\begin{proof}
For $1\le i\le n$, let $W_i\subset\A^n$ be the $\fq$--hypersurface
defined by $X_i^q-X_i$. Then we have $V\cap\fq^n=V\cap
W_1\cap\cdots\cap W_n$. Therefore, applying \cite[Proposition
2.3]{HeSc82} we obtain the inequality $\#(V\cap W_1\cap\cdots\cap
W_n)=\deg(V\cap W_1\cap\cdots\cap W_n)\le \delta q^r$, which
finishes the proof.\end{proof}\newline
We observe that when $r=n-1$, i.e. when $V$ is a hypersurface
defined by a polynomial $f\in\fq[\xon]$, the lemma implies that
the number of $q$--rational zeros of $f$ is at most $dq^{n-1}$.
\begin{lemma}\label{lemma:hesch3} Let $f_1\klk f_s\in\fq[X_1\klk
X_n]$ $(s\ge 2)$ be nonzero polynomials of degree at most
$\delta>0$ without a common factor in $\cfq[\xon]$, and let
$V\subset\A^n$ be the $\fq$--variety defined by $f_1\klk f_s$.
Then $\#(V\cap\fq^n)\le \delta^2 q^{n-2}$.
\end{lemma}
\begin{proof}
Since $f_1,f_2$ have no common factors in $\cfq[X_1\klk\!X_n]$, we
have that $V(f_1,f_2)$ is an $\fq$--variety of dimension $n-2$.
From the B\'ezout inequality (\ref{equation:Bezout}) we conclude
that $\deg V(f_1,f_2)\le \delta^2$ holds. Then Lemma
\ref{lemma:hesch2} shows that $\#\big(V(f_1,f_2)\cap\fq^n\big)\le
\delta^2q^{n-2}$ holds. This implies $\#(V\cap\fq^n)\le
\delta^2q^{n-2}$.\end{proof}\newline
Let us remark that the upper bound of Lemma \ref{lemma:hesch3}
improves the upper bounds $2n\delta^3q^{n-2}$ of \cite[Lemma
4]{Schmidt74} and $\delta^3q^{n-2}$ of \cite[Lemma
IV.3D]{Schmidt76}.
\begin{lemma}\label{lemma:pnai}
Let $V\subset\A^n$ be an $\fq$--irreducible variety of dimension
$r\ge 0$ and degree $\delta$ which is not absolutely irreducible.
Then the inequality $\#(V\cap \fq^n)\leq\delta^2q^{r-1}/4$ holds.
\end{lemma}
\begin{proof} Let $V=V_1\cup\cdots\cup V_s$ be the
decomposition of $V$ into $\cfq$--irreducible components and let
$\delta_i$ denote the degree of $V_i$ for $i=1,\dots,s$. Our
hypotheses imply $s\ge 2$. Since every $q$--rational point of $V$
belongs to $V_i$ for $1\le i\le n$, we see that $V\cap\fq^n\subset
V_1\cap V_2\cap \fq^n$ holds. Therefore, applying Lemma
\ref{lemma:hesch2} we have $\# (V_1\cap V_2\cap \fq^n)\le
\delta_1\delta_2q^{r-1}\le\delta^2q^{r-1}/4$.
\end{proof}\newline
If $V$ is an $\fq$--hypersurface irreducible but not absolutely
irreducible, our estimate gives $\#(V\cap
\fq^n)\leq\delta^2q^{n-2}/4$, which improves the upper bound
$\delta q^{n-1}-(\delta-1)q^{n-2}$ of \cite[Theorem 3.1]{ChRo96},
obtained assuming that $1<\delta<q-1$ holds. Indeed, our upper
bound is valid without any restriction on $q$, and for $\delta\le
q$ we have $\delta^2q^{n-2}/4<\delta q^{n-1}-(\delta-1)q^{n-2}$.
%
%
\section{On the effective first Bertini Theorem.}
\label{section:algoritmoygeneric}
Let be given an absolutely irreducible polynomial $f\in \fq[\xon]$
of degree $\delta>0$ and let $H\subset\A^n$ be the affine
$\fq$--hypersurface defined by $f$. Our estimates on the number of
$q$--rational points of $H$ rely on an analysis of the varieties
obtained by intersecting $H$ with an affine linear $\fq$--variety
of dimension 2 ($\fq$--plane for short). For this purpose, we need
an estimate on the number of $\fq$--planes $L$ for which $H\cap L$
has an absolutely irreducible $\fq$--component of degree at most
$D$, for a given $1\leq D\leq \delta-1$.

Following \cite{Kaltofen95a}, we analyze the genericity condition
underlying the nonexistence of irreducible components of $H\cap L$
of degree at most $D$. In order to do this, in the next section we
introduce an algorithm which, given a bivariate polynomial $f \in
K[X,Y]$,  finds the irreducible factors of $f$ over $K$  of degree
at most $D$. Then, in Section \ref{section:genericity condition}
we obtain a suitable upper bound on the genericity condition we
are considering.
%
%
\subsection{An algorithm computing the irreducible factors  of degree at most
$D$ of a bivariate polynomial over a field
$K$.}\label{section:algoritmo}
The algorithm we exhibit in this section is a variant of the
corresponding algorithm of \cite{Kaltofen95a}.

{\bf Algorithm} {\em Factorization over the Coefficient Field of
degree at most $D$.}

{\em Input:} A polynomial $f\in K[X,Y]$ monic in $X$ of degree at
most $\delta$, where $K$ is an arbitrary field, such that the
resultant $Res_X (f(X,0),\partial f(X,0)/\partial X)\neq 0$, and
an integer $D$ with $1\le D\le \delta-1$.

{\em Output:} Either the  algorithm returns the list of
irreducible factors of $f$ defined over $K$ of degree at most $D$,
or $f$ will not have irreducible factors in $K[X,Y]$ of degree at
most $D$.

{\em Set the maximum order of approximation needed:
$\ell_{max}\leftarrow 2D\delta.$ \medskip\newline
For all roots $\zeta_i\in \overline{K}$ of $f(X,0)\in K[X]$ Do
steps {\bf N} and {\bf L}}. \medskip\newline
{\bf Step N:} Let $K_i:=K(\zeta_i)$ and {\em Set the initial
points for the Newton iteration}
$$\alpha_{i,0}\leftarrow \zeta_i \in K_i,\ \beta_{i,0}\leftarrow
{(\partial f/\partial X)(\alpha_{i,0},0)}^{-1} \in K_i.$$
${}^{}$\quad(Now we perform Newton iteration) \medskip\newline
${}^{}$\quad {\em For $j \leftarrow 0,\ldots,\lfloor
\log_2(\ell_{\max})\rfloor$ Do}

\vskip-20pt
$$\begin{array}{l}\alpha_{i,j+1}\leftarrow(\alpha_{i,j}-
\beta_{i,j}f(\alpha_{i,j},Y))\!\!\!\!\pmod {Y^{2^{j+1}}}\\
\beta_{i,j+1}\leftarrow \left(2\beta_{i,j}-(\partial f/\partial
X)(\alpha_{i,j+1},Y) \beta_{i,j}^2 \right)\!\!\!\!\pmod
{Y^{2^{j+1}}}.\end{array}$$

${}^{}$\quad(Observe that $\alpha_{i,j+1},\beta_{i,j+1}$ are
polynomials of $K_i[Y]$ satisfying \newline ${}^{}\quad
f(\alpha_{i,j+1},Y)\equiv 0 \!\!\!\!\pmod {Y^{2^{j+1}}},
\beta_{i,j+1}\cdot(\partial f/\partial x)(\alpha_{i,j+1},Y) \equiv
1\!\!\!\!\pmod {Y^{2^{j+1}}}.$)\medskip\newline
${}^{}$\quad {\em Set the approximate root} $$\alpha_i\leftarrow
\alpha_{i,\lfloor \log_2(\ell_{\max})\rfloor +1}\!\!\!\!\pmod
{Y^{\ell_{\max}+1}} \in K_i[Y].$$
${}^{}$\quad (Next, we compute the powers of $\alpha_i$.)
\medskip\newline
${}^{}$\quad {\em For $\mu\leftarrow 0,\ldots,\delta-1$ Do}
$$\sum_{k=0}^{\ell_{\max}}a_{i,k}^{(\mu)}Y^k \leftarrow
\alpha_i^{\mu}\!\!\!\!\pmod{Y^{\ell_{\max}+1}}\ \mathrm{with}\
a_{i,k}^{(\mu)} \in K_i.$$

{\bf Step L:} We find the lowest degree polynomial in $K[X,Y]$
whose root is $\alpha_i$.

\vskip-7pt ${}^{}$\quad {\em For $m\leftarrow 1,\ldots,D$ Do}
\medskip\newline
${}^{}$\quad {\em We fix the  order of approximation}:
$\ell\leftarrow 2m\delta.$
\medskip\newline
${}^{}$\quad (We examine if the equation $\alpha_i^m +
\sum_{\mu=0}^{m-1} h_{i,\mu}(Y)\alpha_i^{\mu}\equiv 0
\!\!\!\!\pmod{Y^{\ell +1}}$ has a \newline ${}^{}$\quad solution
for $h_{i,\mu}(Y) \in K[Y]$ with $\deg (h_{i,\mu}) \leq m-\mu $.
Writing $h_{i,\mu}(Y)=$
\linebreak ${}^{}$\quad
$\sum_{\eta=0}^{m-\mu}u_{i,\mu,\eta}Y^{\eta}$, with
$u_{i,\mu,\eta} \in K$, and  collecting the coefficients of $Y^k$
we are \linebreak ${}^{}$\quad led to the following problem.)
\medskip\newline${}^{}$\quad {\em For $0\leq k\leq \ell$, solve
the following linear system over $K$ in the variables\linebreak
${}^{}\quad u_{i,\mu,\eta}$ $(0\leq \mu \leq m-1$, $0\leq \eta
\leq m-\mu)$:}
\begin{equation} \label{sistema}
a_{i,k}^{(m)}+ \sum_{\mu=0}^{m-1} \sum_{\eta=0}^{m-\mu}a_{i,k-\eta
}^{(\mu)}u_{i,\mu,\eta}=0 \ \ (where\ a_{i,\nu}^{(\mu)}=0\ for\
\nu <0),
\end{equation}
\vskip-20pt
%
${}^{}$\quad (Since $\deg h_{i,\mu}\leq m-\mu$ holds, for every
$\mu$ we have $m-\mu+1$ indeterminates, \linebreak ${}^{}$\quad
which implies that the system has $(m+1)(m+2)/2-1$
indeterminates.)
\medskip\newline${}^{}$\quad {\em If $(\ref{sistema})$
has a solution then}
$$f_i(X,Y)\leftarrow X^m + \sum_{\mu=0}^{m-1}
\sum_{\eta=0}^{m-\mu}u_{i,\mu,\eta}Y^{\eta}X^{\mu}.$$
\medskip\newline${}^{}$\quad (The polynomial $f_i(X,Y)$ is an
irreducible factor of $f(X,Y)$ of degree $D$ or \linebreak
${}^{}$\quad some factor of it is an irreducible factor of degree
less than $D$.)
\medskip\newline${}^{}$\quad {\em Check if $f_i$ has
been produced by a root $\zeta_l$ with $l<i$. If not, add $f_i$ to
the list\linebreak ${}^{}$\quad of irreducible factors of degree
less than $D$.}
\medskip\newline${}^{}$\quad {\em If $(\ref{sistema})$ has no
solution for all $i=1,\ldots,\delta$ and $m=1\klk D$, then $f$ has
no \linebreak ${}^{}$\quad irreducible factors in $K[X,Y]$ of
degree at most $D$.}

The next lemma proves the correctness of this algorithm:
\begin{lemma}
The  polynomial $f(X,Y)$ has  an irreducible factor over $K$ of
degree at  most $D$ if and only if at least one of the $D\delta$
linear systems (\ref{sistema}) has a  solution in  $K$.
\end{lemma}
\begin{proof}
Suppose that (\ref{sistema}) has a  solution in $K$, i.e. there
exists $1\le i\le \delta$ and a polynomial $g_i(X,Y)\in K[X,Y]$ of
degree $1\le m\le D$ such that $g(\alpha_i,Y)\equiv 0 \
(\mathrm{mod}\ Y^{2m\delta+1})$. Let $\rho\in K[Y]$ denote the
resultant $\rho(Y):=\mbox{Res}_X(f,g)$. Evaluating $\rho$ at
$X=\alpha_i$ we conclude that $\rho(Y)\equiv 0 \ (\mathrm{mod}\
Y^{2m\delta +1})$. Since $\rho$ is a polynomial of degree at most
$2m\delta$, we conclude that $\rho=0$ holds. Hence $\gcd(f,g)$ is
a non--trivial element of $K[X,Y]$, and therefore is a factor of
$f$ of degree at most $D$.

Now, suppose that $f(X,Y)$ has an irreducible factor $g(X,Y)\in
K[X,Y]$ of degree at most $D\geq 1$. Then there exists a
non--trivial factorization $f(X,Y)=g(X,Y)h(X,Y)$ over $K[X,Y]$.
Let $1\le i\le d$ be an integer for which $g(\alpha_{i,0},0)=0$.
Then $h(\alpha_{i,0},0)\not= 0$, which implies
$h(\alpha_{i,0},Y)\not\equiv 0\ (\mathrm{mod}\ Y)$ and thus
$h(\alpha_{i,j},Y)\not\equiv 0\ (\mathrm{mod}\ Y)$ for $1\le j\le
\lfloor\log_2(\ell_{\mathrm{max}})\rfloor+1$. Therefore, we have
$h(\alpha_i,Y)\equiv 0\ (\mathrm{mod}\ Y)$, which combined with
$f(\alpha_i,Y)\equiv 0 \ (\mathrm{mod}\ Y^{2D\delta+1})$ shows
that $g(\alpha_i,Y)\equiv 0 \ (\mathrm{mod}\ Y^{2D\delta+1})$
holds. We conclude that the coefficients of $g$, considered as
polynomial of $K[Y][X]$, furnish a solution to at least one of the
$D\delta$ linear systems (\ref{sistema}). This completes the
proof.
\end{proof}
%
%
\subsection{The genericity condition underlying the existence of
irreducible components of a given degree.}
\label{section:genericity condition}
The estimates on the number of $q$--rational points of a given
absolutely irreducible $(n-r)$--dimensional $\fq$--variety
$V\subset\A^n$ of e.g. \cite{Schmidt76}, \cite{HuWo98},
\cite{CaMa02} depend strongly on a suitable effective version of
the first Bertini Theorem. As it is well--known, the first Bertini
Theorem (see e.g. \cite[\S II.6.1, Theorem 1]{Shafarevich94}
asserts that the intersection $V\cap L$ of $V$ with a generic
affine variety variety $L\subset\A^n$ of dimension $r+1$ is an
absolutely irreducible curve. An effective version of the first
Bertini Theorem aim at estimating the number of planes $L$ for
which $V\cap L$ is not an absolutely irreducible curve, and is
usually achieved by analyzing the genericity condition underlying
the choice of a linear variety $L$. The estimates for
hypersurfaces we shall present in the next sections rely on a
variant of the effective first Bertini Theorem, which estimates
the number of planes $L$ whose intersection with a given
absolutely irreducible $\fq$--hypersurface $H$ has absolutely
irreducible $\fq$--components of degree at most $D$ for a given
$1\le D\le \delta-1$.

For this purpose, let $f \in K[\xon]$ be an  absolutely
irreducible polynomial of degree $d$ and let be given $1\le D\le
\delta-1$. For $\nu_1\klk \nu_n, \omega_2\klk \omega_n \in
\overline{K}$, we consider the polynomial $$\chi(X,Y,Z_2\klk
Z_n):= f(X+\nu_1,\omega_2X + Z_2 Y +\nu_2\klk \omega_nx + Z_n Y
+\nu_n)$$ as an element of $\overline{K}[X,Y,Z_2\klk Z_n]$.
Following \cite[Lemmas 4 and 5]{Kaltofen95a}, there exists a
non--zero polynomial $\Upsilon \in K[V_1\klk V_n, W_2\klk W_n]$ of
degree at most $2\delta^2$ such that for any
$\nu_1,\ldots,\nu_n,\omega_2,\ldots,\omega_n \in \overline{K}$
with
\begin{equation}\label{Upsilon}
\Upsilon(\nu_1\klk \nu_n, \omega_2\klk \omega_n)\neq 0
\end{equation}
 the following conditions are satisfied:
\begin{itemize}
\item the leading coefficient of $\chi$ with respect to $X$ is a
non--zero element of $\overline{K}$,
\item the discriminant of $\chi(X,0,Z_2\klk Z_n)$  with respect to $X$  is
non--zero,
\item $\chi$ is an irreducible element of $\overline{K}[X,Y,Z_2\klk Z_n]$.
\end{itemize}

Under the assumption of condition (\ref{Upsilon}), Kaltofen proves
a crucial fact for his effective version of the first Bertini
Theorem \cite[Theorem 5]{Kaltofen95a}: he shows the existence of a
polynomial $\Psi \in \overline{K}[Z_2\klk Z_n]$ of degree at most
$3\delta^4/2 -2 \delta^3 +\delta^2/2$ such that for any
$\eta:=(\eta_2\klk \eta_n) \in \overline{K}^{n-1}$ with
$\Psi(\eta)\not=0$, the bivariate polynomial $\chi(X,Y,\eta_2\klk
\eta_n)\in \overline{K}[X,Y]$ is absolutely irreducible.

Hence, for $\Xi:=\Upsilon(V_1\klk V_n, W_2\klk W_n)\Psi(Z_2\klk
Z_n)$ we have $\deg \Xi\le 3\delta^4/2-2\delta^3+ 5\delta^2/2$ and
for any
$(\nu,\omega,\eta):=(\nu_1,\ldots,\nu_n,\omega_2,\ldots,\omega_n,
\eta_2\klk \eta_n) \in \overline{K}^{3n-2}$ with
$\Xi(\nu,\omega,\eta)\not=0$, the polynomial $\chi(X,Y,\eta_2\klk
\eta_n):= f(X+\nu_1,\omega_2X+\eta_2Y+\nu_2\klk \omega_nX+\eta_n
Y+\nu_n)$ is absolutely irreducible. In particular, for $K=\fq$ we
deduce the following corollary:
\begin{corollary}\label{corokaltofen}
Let $f \in \fq[\xon]$ be absolutely irreducible of degree
$\delta>0$. Then there exists at most
$(3\delta^4/2-2\delta^3+5\delta^2/2)q^{3n-3}$ elements
$(\nu,\omega,\eta)\in \fq^{3n-2}$ for which $\chi(X,Y,\eta)$ is
not absolutely irreducible.
\end{corollary}
Our goal is to obtain a degree estimate, similar to that of $\deg
\Xi$, on the genericity condition underlying the nonexistence of
absolutely irreducible factors of $\chi(X,Y,\eta_2\klk \eta_n)$ of
degree at most $D$ for a given $1\leq D\leq \delta-1$. Our next
theorem, a variant of \cite[Theorem 5]{Kaltofen95a}, will be a
crucial point for our estimates of the following sections.
\begin{theorem}\label{nuestrokaltofen}
Let  $1\leq D \leq \delta-1$ and suppose that
$\nu_1,\ldots,\nu_n,\omega_2,\ldots,\omega_n$ satisfy condition
(\ref{Upsilon}). Then there exists a nonzero polynomial $\Psi_D
\in \overline{K}[Z_2,\ldots,Z_n]$ of degree $$\deg \Psi_D \le
D\delta^2(D+1)(D+2)-(D^2+3D)(D^2+3D+2)\delta/8$$ such that for any
$\eta:=(\eta_2\klk\eta_n)\in \overline{K}^{n-1}$ with
$\Psi_D(\eta)\neq 0$, the polynomial $\chi(X,Y,\eta):=
f(X+\nu_1,\omega_2X+\eta_2Y+ \nu_2,\ldots,
\omega_nX+\eta_nY+\nu_n)$ has no irreducible factors of degree at
most $D$ in $\overline{K}[X,Y]$.
\end{theorem}
\begin{proof} Since by assumption $\chi$ is irreducible over
$\overline{K}[Z_2,\ldots,Z_n][X,Y]$, Gauss Lemma implies that
$\chi$ is irreducible over $\overline{K}(Z_2,\ldots,Z_n) [X,Y]$.
Therefore, applying  algorithm {\em Factorization over the
Coefficient Field of degree at most $D$} to the polynomial
$\psi:={l}^{-1}\chi \in \overline{K}(Z_2,\ldots,Z_n)[X,Y],$ where
$l\in\overline{K}$ is the leading coefficient of $\chi$, since
$\psi(X,0) \in \overline{K}[X]$, the root $\zeta_i$ used to
construct the field $K_i:=K(\zeta_i)$ of the algorithm is actually
an element of $\overline{K}$ for $1\le i\le \delta$. Then the
irreducibility of $\psi$ over $\overline{K}(Z_2,\ldots,Z_n)[X,Y]$
implies that the linear system (\ref{sistema}) derived in step L
has no solution in the field $K_i$ for $1\le m\le D$ and $1\le
i\le \delta$. This implies that for $m=D$ and $1\le i\le \delta$,
the augmented matrix of the system,
$\widetilde{M}^{(i)}_D(Z_2,\ldots,Z_n)$, has rank greater than
that of the matrix of the coefficients
$M^{(i)}_D(Z_2,\ldots,Z_n)$. Since $\partial{\psi}(X,0)/\partial X
\in \overline{K}$, all the denominators used in the  construction
of this system are elements of $\overline{K}$. Let $\Psi^{(i)}_D
\in \overline{K}[Z_2,\ldots,Z_n]$ be a maximal nonzero minor of
the augmented matrix $\widetilde{M}^{(i)}_D(Z_2,\ldots,Z_n)$ and
let $\eta:=(\eta_2,\ldots,\eta_n)\in \overline{K}^{n-1}$ satisfy
$\prod_{i=1}^{\delta}\Psi^{(i)}_D(\eta)\neq 0$. Then the
specialized system (\ref{sistema}) has no solutions for $i=1\klk
\delta$, which implies that $\chi(X,Y,\eta_2,\ldots,\eta_n)$ has
no irreducible factors over $\overline{K}[X,Y]$ of degree at most
$D$, because algorithm {\em Factorization over the Coefficient
Field of degree at most $D$} fails to find such a factor of
$\chi(X,Y,\eta)$ over $\overline{K}$. Therefore, $\Psi_D:=
\prod_{i=1}^{\delta}\Psi^i_D$ is the polynomial we are looking
for.

Now we show that the degree estimate for $\Psi_D$ holds. The
degree estimate essentially follows from the degree estimate of
the proof of \cite[Theorem 5]{Kaltofen95a}, taking into account
that we have a different number of indeterminates and a different
order of approximation. Indeed, for every root $\zeta_i$ of
$\psi(x,0)$, the corresponding linear system for $m=D$ has
$(D+1)(D+2)/2-1$ indeterminates. Hence, any maximal nonzero minor
$\Psi^{(i)}_D$ satisfies the following degree estimate:
\[\begin{array}{ccl}
  \deg_{Z_2,\ldots,Z_n}\!\Psi^{(i)}_D & \le &
   \sum_{j=0}^{(D+1)(D+2)/2-1}(\ell_{\max}-j) \\
   & \leq &\displaystyle 2D\delta(D+1)(D+2)/2-(D^2+3D)(D^2+3D+2)/8 \\
   & = &\displaystyle D(D^2+3D+2)(\delta-(D+3)/8).
\end{array}\]
This immediately implies the degree estimate for $\Psi_D$ of the
theorem.\end{proof}

Since Theorem \ref{nuestrokaltofen} is valid under the assumption
of condition (\ref{Upsilon}), if we define
$$\Xi_D:=\Upsilon(V_1,\ldots,V_n,W_2,\ldots,W_n)
\Psi_D(Z_2,\ldots,Z_n),$$
then $\Xi_D$ is a polynomial in $3n-2$ indeterminates with
coefficients in $\overline{K}$ of degree bounded by
$$\deg \Xi_D \le D^3\delta^2- D^4\delta/8-3D^3\delta/4 +
3D^2\delta^2 -11D^2\delta/8 +2D\delta^2 -3D\delta/4 + 2\delta^2,$$
which satisfies the following property: for every
$(\nu,\omega,\eta) \in \overline{K}^{3n-2}$ with \linebreak $\Xi_D
(\nu,\omega,\eta)\neq 0$, the polynomial
$\chi(x,y,\eta_2,\ldots,\eta_n)$ has no irreducible factors over
$\overline{K}$ of degree at most $D$. Therefore, for $K=\fq$ we
have the following corollary:
\begin{corollary}\label{corokaltofenD}
Let $f \in \fq[\xon]$ be an absolutely irreducible polynomial of
degree $\delta\geq 2$ and let be given an integer $D$ with $1\leq
D\leq \delta-1$. Then there are at most
$(D^3\delta^2-D^4d/8-{3}D^3\delta/{4}+3D^2\delta^2-{11}D^2\delta/{8}+2D\delta^2
-{3}D\delta/{4}+2\delta^2)q^{3n-3}$ elements $(\nu,\omega,\eta)\in
\fq^{3n-2}$ for which $\chi(X,Y,\eta_2\klk \eta_n)$ has a
nonconstant irreducible factor over $\cfq[X,Y]$ of degree at most
$D$.
\end{corollary}
%
%
\section{On the intersection of an absolutely irreducible
$\fq$--hypersurface with an $\fq$--plane.} \label{section:number
of planes}
Following \cite{Schmidt74}, in this section we estimate the number
of planes whose restriction to a given polynomial $f$ has a fixed
number of $\fq$--absolutely irreducible factors. For this purpose,
we shall consider the results from the previous section from a
geometric point of view. We shall work with the field $K:= \fq$
and every non--zero $(\nu,\omega,\eta)\in \fq^{3n-2}$ shall be
considered as providing a parametrization of a linear affine
$\fq$--variety of $\A^n$ of dimension 2.

More precisely, let be given a polynomial $f \in \fq[\xon]$ of
degree $\delta >0$. For an affine linear $\fq$--variety $L\subset
\A^n$ of dimension 2 (an $\fq$--plane for short), we represent the
restriction of $f$ to $L$ as a bivariate polynomial $f_L \in
\fq[X,Y]$, where $X,Y$ are the parameters of a given
parametrization of $L$. Let us remark that for every such $L$, the
polynomial $f_L$ is univocally defined up to an $\fq$--definable
affine linear change of coordinates. Therefore, its degree and
number of absolutely irreducible components do not depend on the
particular parametrization of $L$ we choose to represent $L$ (cf.
\cite[V.\S 4]{Schmidt76}).

In particular, we shall be concerned with the $\fq$--planes of
$\A^n$ which have an $\fq$--definable parametrization of the
following type:
\begin{equation}\label{equation:paramet}
\begin{array}{rc}
 X_1= \nu_1\! + X,& \quad \quad X_i=\nu_i+\omega_iX+\eta_iY\quad (2\le i\le n).
\end{array}\end{equation}
Let $M_T^{(2)}$ denote the set of all $\fq$--planes of $\A^n$ and
let $M^{(2)}$ denote the subset formed by the elements of
$M_T^{(2)}$ having a parametrization as in
(\ref{equation:paramet}).

Our purpose is to analyze the number of  absolutely irreducible
factors of $f_L$ for a given $\fq$--plane $L \subset \A^n$. Hence,
for a plane $L \in M^{(2)}$, let  the non--negative integer
$\nu(L)$ represent the number of absolutely irreducible
$\fq$--definable factors of the  polynomial $f_L$. Here $0\leq
\nu(L)\leq \delta $ if $f$ does not vanish identically on $L$ and
$\nu(L)=q$ otherwise.  Further, let $\Pi_{j}$ be the set of planes
$L$ with $|\nu(L)-1|=j$. Thus, $\Pi_1$ is the set of planes $L$
with 0 or 2 absolutely irreducible $\fq$--factors, $\Pi_j$ is the
set of planes $L$ for which $f_L$ has $j+1$ absolutely irreducible
$\fq$--factors for $j=0,2\klk \delta-1$, and $\Pi_{q-1} $ is the
set of planes $L$ for which $f_L$ vanish identically.

Observe that if $L\in \Pi_j$ for a given $j=0\klk \delta-1$ then
$f_L$ has an absolutely irreducible factor of degree at most
$D_j:=\lfloor{\delta}/{(j+1)}\rfloor$. Theorem
\ref{nuestrokaltofen} asserts that for any plane $L \in M^{(2)}$
with an $\fq$--parametrization as in (\ref{equation:paramet}) for
which $\Xi_{D_j}(\nu,\omega,\eta) \neq 0$ holds, the polynomial
$f_L$ has no irreducible factors over $\cfq$ of degree at most
$D_{j}$.
Therefore, for every such $(\nu,\omega,\eta)\in\fq^{3n-2}$ the
polynomial $f_L$ cannot have at least $j+1$ irreducible factors
over $\cfq$, which in particular implies that $L \notin
\Pi_{j}\cup \cdots \cup \Pi_{\delta-1}$ holds. Hence we have
$$\Pi_{j}\cup \cdots \cup \Pi_{\delta-1} \subset
\{(\nu,\omega,\eta)\in \fq^{3n-2}:
\Xi_{D_j}(\nu,\omega,\eta)=0\}.$$

Taking into account that every plane of $M^{(2)}$ has $q^3(q-1)$
$\fq$--parametrizations as in (\ref{equation:paramet}), from Lemma
\ref{lemma:hesch2} we deduce the following estimate:
$$\#(\Pi_{j})+ \cdots + \#(\Pi_{\delta-1}) \leq \deg
\Xi_{D_j}\frac{q^{3n-3}}{q^3(q-1)}.$$

Therefore, from Corollary \ref{corokaltofenD} we obtain the
estimate
\begin{equation}\label{sumadepi}\begin{array}{c}\!\!\!\!\!\!
\!\!\!\!\!\! {\displaystyle
\sum_{k=j}^{\delta-1}}\#(\Pi_{k})\!\!\le \!
\!\Big(\delta^5\Big(\frac{1}{j^3}-
\frac{1}{8j^4}\Big)\!+\!3\delta^4\Big(\frac{1}{j^2}-\frac{1}{4j^3}\Big)
\!+\!\delta^3\Big(\frac{2}{j}-\frac{11}{8j^2}\Big)
-\frac{3}{4}\frac{\delta^2}{j} \!+ \!2\delta^2
\Big)\frac{q^{3n-6}}{(q-1)}.\end{array}\end{equation}
The following proposition is crucial for our estimates for an
absolutely irreducible $\fq$--hypersurface of the next section. It
yields a better estimate than that obtained by a straightforward
application of Theorem \ref{corokaltofen}.
\begin{proposition}\label{prop:cota}
Let $f \in \fq[\xon]$ be a polynomial of degree $\delta >1 $. Then
the following estimate holds:
$$\sum_{j=1}^{\delta-1}j\#(\Pi_j) \leq (2 \delta^{\frac{13}{3}}+
3\delta^{\frac{11}{3}})\frac{ q^{3n-3}}{q^3(q-1)}.$$
\end{proposition}
\begin{proof}
For $\delta=2$ the expression $\sum_{j=1}^{\delta-1}j\#(\Pi_j)$
consists of only one term, namely $\Pi_1$, and therefore Corollary
\ref{corokaltofen} yields
\begin{equation}\label{d=2}\#(\Pi_1)\leq \left(\frac{3}{2}\delta^4 - 2\delta^3 +
\frac{5}{2}\delta^2\right)\frac{q^{3n-6}}{(q-1)}\leq
\frac{3}{2}\delta^4\frac{q^{3n-6}}{(q-1)}.\end{equation}

Hence, we may assume without loss of generality $\delta\geq3$. Let
$r$ be a real number, to be fixed below, lying in the open
interval $(1,\delta-1)$. We have:
$$\sum_{j=1}^{\delta-1}j \# (\Pi_j)\!= \!\sum_{j=1}^{\lfloor r
\rfloor}j\#(\Pi_j)+ \sum_{j=\lfloor r \rfloor +1}^{\delta-1}\!\!j
\#(\Pi_j)\! \leq \!\lfloor r \rfloor
\sum_{j=1}^{\delta-1}\#(\Pi_j) + \sum_{j=\lfloor r \rfloor
+1}^{\delta-1}\!\!(j-\lfloor r \rfloor)\#(\Pi_j).$$
By Corollary \ref{corokaltofen} we have:
\begin{equation}\label{1termino}\lfloor r \rfloor
\sum_{j=1}^{\delta-1}\#\Pi_j \leq r\left(\frac{3}{2}\delta^4 -
2\delta^3 +
\frac{5}{2}\delta^2\right)\frac{q^{3n-6}}{(q-1)}.\end{equation}
On the other hand, by inequality (\ref{sumadepi}) we have:
\begin{equation}\label{2termino}\begin{array}{lll}
\displaystyle \sum_{j=\lfloor r \rfloor +1}^{\delta-1}(j-\lfloor r
\rfloor)\#(\Pi_j)& = & \displaystyle \sum_{j=\lfloor r \rfloor
+1}^{\delta-1}(\#(\Pi_j)+\cdots+\#(\Pi_{\delta-1}))\\ &\leq &
\left( \delta^5 c_1 +3\delta^4c_2  + \delta^3c_3 -
\frac{3}{4}\delta^2c_4 + 2\delta^2 \right)\!\!\ds \frac{
q^{3n-6}}{(q-1)},
\end{array}\end{equation}
where $c_1,c_2,c_3,c_4$ are the following numbers:
$$ c_1:=\!\!\!\sum_{j=\lfloor r \rfloor
+1}^{\delta-1}\frac{1}{j^3}- \frac{1}{8j^4},\,
c_2:=\ds\!\!\!\sum_{j=\lfloor r \rfloor
+1}^{\delta-1}\frac{1}{j^2}-\frac{1}{4j^3},\,
c_3:=\!\!\!\sum_{j=\lfloor
r\rfloor+1}^{\delta-1}\frac{2}{j}-\frac{11}{8 j^2},\,
c_4:=\!\!\!\sum_{j=\lfloor r\rfloor+1}^{\delta-1}\frac{1}{j}.$$

We observe that any decreasing positive real function $g$
satisfies the inequality $\sum_{j=\lfloor r \rfloor
+1}^{\delta-1}g(j) \leq \int_{ r}^{\delta-1}g(x)dx$. Let
$r:=\delta^{\frac{1}{3}}$. Using the fact that
$1<\delta/(\delta-1)\leq 3/2$ holds for $\delta\geq 3$, we have
the following inequalities:
$$\begin{array}{rcl} \delta^5c_1& \leq & \delta^5  \left
(\frac{1}{2}(\delta^{\frac{1}{3}})^{-2} - \frac{1}{24}
(\delta^{\frac{1}{3}})^{-3}-\frac{1}{2}(\delta-1)^{-2}
+\frac{1}{24}(\delta-1)^{-3}\right) \\ & \leq &
\frac{1}{2}\delta^{\frac{13}{3}} - \frac{1}{24}\delta^4
-\frac{1}{2}\delta^3 + \frac{9}{64}\delta^2 \\ & \\3\delta^4c_2&
\leq & 3\delta^4\left(\delta^{-\frac{1}{3}} -
\frac{1}{8}(\delta^{\frac{1}{3}})^{-2} -(\delta-1)^{-1}+
\frac{1}{8}(\delta-1)^{-2} \right)\\& \leq & 3
\delta^{\frac{11}{3}}-\frac{3}{8}\delta^{\frac{10}{3}} -3\delta^3
+ \frac{27}{32}
 \delta^2\\& \\ \delta^3c_3& \leq &\delta^3(2\ln(\delta-1)+ \frac{11}{8}(\delta-1)^{-1} -2\ln
\delta^{\frac{1}{3}} - \frac{11}{8}\delta^{-\frac{1}{3}}) \\ &\leq
& \frac{4}{3}\delta^3\ln \delta + \frac{33}{16} \delta^2 -
\frac{11}{8}\delta^{\frac{8}{3}}. \end{array} $$

This implies that the following estimate holds:
\begin{equation}\label{3termino}\begin{array}{l}
\delta^5c_1 +3\delta^4c_2 +\delta^3c_3 - \frac{3}{4}\delta^2c_4
+2\delta^2\le\\ \le\frac{1}{2}\delta^{\frac{13}{3}} -
\frac{1}{24}\delta^4 +3\delta^{\frac{11}{3}}
-\frac{3}{8}\delta^{\frac{10}{3}} +\frac{4}{3}\delta^3\ln \delta
-\frac{7}{2}\delta^3 -\frac{11}{8}\delta^{\frac{8}{3}} +
\frac{323}{64}\delta^2
\end{array}
\end{equation}

Now, putting together (\ref{1termino}), (\ref{2termino}) and
(\ref{3termino}), and taking into account that
$\frac{4}{3}\delta^3\ln \delta\le 3\delta^{3+\, 1/3}$ holds for
$\delta\geq 3$, we obtain:
$$\begin{array}{rcl}\sum_{j=1}^{\delta-1}j\#(\Pi_j)&\le& 2\delta^{
\frac{13}{3}} - \frac{1}{24}\delta^4 +3d^{\frac{11}{3}}
+\frac{5}{8}\delta^{\frac{10}{3}} -\frac{7}{2}\delta^3
-\frac{11}{8}\delta^{ \frac{8}{3}}
+\frac{5}{2}\delta^{\frac{7}{3}} +\frac{323}{64}\delta^2 \\ &\le &
2 \delta^{\frac{13}{3}}+ 3 \delta^{\frac{11}{3}} \end{array}$$ for
$\delta\geq 3$. This proves the proposition.
\end{proof}
%
%
\section{Estimates for an absolutely irreducible
$\fq$--hypersurface.}\label{section:hypersurface}
In this section we obtain different types of estimates on the
number of $q$--rational points of a given absolutely irreducible
$\fq$--hypersurface. First we exhibit an estimate which holds
without any regularity condition and  improves
(\ref{Ghorpadehypersurface}) and (\ref{estimate:HuangWong}). Then
we show an estimate which improves both the right--hand side and
the regularity condition of the lower bound
(\ref{estimate:Schmidt74}), providing also a corresponding upper
bound. Finally, we extend these estimates to the case of an
arbitrary $\fq$--hypersurface.

For this purpose, we shall follow an approach that combines both
ideas of \cite{Schmidt74} and \cite{Schmidt76} with the estimate
of the previous section. This approach is based on estimating the
number of $q$--rational points lying in the intersection of a
given absolutely irreducible $\fq$--hypersurface with all the
$\fq$--planes of $\A^n$.
%
%
\subsection{An estimate without any regularity condition.}
\label{subsec:hypersup_without_regularity}
In what follows, we shall apply the following lemma from
\cite{Schmidt74}.
\begin{lemma}{\rm \cite[Lemma 5]{Schmidt74}}\label{lemaSchmidt}
Let  $f\in \fq[X,Y]$ be a polynomial of degree $\delta>0$ and let
$\nu$ be the number of distinct absolutely irreducible
$\fq$--definable factors of $f$. Then the number $N$ of zeros of
$f$ in $\fq^2$ satisfies $$|N - \nu q| \leq \omega(q,\delta) +
\delta^2, $$ where $\omega(q,\delta):= (\delta-1)(\delta-2)q^{1/2}
+ \delta+1$.
\end{lemma}

Let be given an absolutely irreducible polynomial $f \in
\fq[\xon]$ of degree $\delta>0$. Recall that $M_T^{(2)}$ and
$M^{(2)}$ denote the number of $\fq$--planes of $\A^n$ and the
number of $\fq$--planes with a parametrization as in
(\ref{equation:paramet}) respectively. Further, for $j=0,2,\dots,
\delta-1$ let $\Pi_j$ be the number of $\fq$--planes $L\in
M^{(2)}$ for which the restriction $f_L$ of $f$ to $L$ has $j+1$
absolutely irreducible $\fq$--definable factors, let $\Pi_1$ be
the number of $\fq$--planes $L\in M^{(2)}$ for which $f_L$ has 0
or 2 absolutely irreducible $\fq$--factors, and let $\Pi_{q-1}$
denote the number of $\fq$--planes $L\in M^{(2)}$ for which $f_L$
vanishes identically. Let us introduce the following quantities:
$$\begin{array}{lccr} A:=\#\,M^{(2)}, & \quad B:=\displaystyle
\sum_{j=1}^{\delta-1}j\#(\Pi_j),& \quad C:=\#(\Pi_{q-1}), & \quad
D:=\# M_T^{(2)}- \# M^{(2)}.\end{array}$$
Let $E$ denote the number of elements of $M_T^{(2)}$ containing a
given point of $\fq^n$.

We recall that any element of $M^{(2)}$ is represented by
$D':=q^3(q-1)$ different parametrizations of type
(\ref{equation:paramet}). Therefore, taking into account that
there are $q^{2n-1}(q^{n-1}-1)$ different parametrizations of type
(\ref{equation:paramet}), we conclude that
\begin{equation}\label{eq:estimate_A}
A=\frac{q^{2n-1}(q^{n-1}-1)}{q^3(q-1)}
\end{equation}
holds. By Proposition \ref{prop:cota} we have
$B\le(2\delta^{\frac{13}{3}}+ 3\delta^{\frac{11}{3}})\frac{
q^{3n-3}}{q^3(q-1)}$, which implies
\begin{equation}\label{eq:estimate_B/A}
\frac{B}{A}\le\left(2\delta^{\frac{13}{3}}+
3\delta^{\frac{11}{3}}\right)\frac{ q^{n-2}}{q^{n-1}-1}.
\end{equation}
By a simple recursive argument we may assume without loss of
generality that $f$ cannot be expressed as a polynomial in $n-2$
variables (see e.g. \cite{Schmidt76}). Let us fix $c\in\fq^{n-2}$
for which $f(c,X_{n-1},X_n)$ vanishes identically. Let us write
$f=\sum_{\alpha\in\mathcal{J}}f_\alpha X_{n-1}^{\alpha_1}
X_n^{\alpha_2}$, where $\mathcal{J}\subset( \Z_{\ge 0})^2$ is a
suitable finite set and $f_\alpha\in\fq[\xo{n-2}]$ for any
$\alpha=(\alpha_1,\alpha_2)\in \mathcal{J}$. Since $f$ is not a
polynomial of $\fq[\xo{n-2}]$, it follows that $f_\alpha(c)=0$ for
any $\alpha\in\mathcal{J}$. By the absolute irreducibility of $f$
we have that the set of polynomials
$\{f_\alpha:\alpha\in\mathcal{J}\}\subset\fq[\xo{n-2}]$ does not
have nontrivial common factors in $\cfq[X_1\klk X_{n-2}]$. Then
Lemma \ref{lemma:hesch3} implies that there exist at most
$\delta^2q^{n-4}$ elements $c\in\fq^{n-2}$ for which
$f(c,X_{n-1},X_n)=0$ holds, and hence there exist at most
$\delta^2q^{n-4}$ linear varieties $L$ of $M_2$ parallel to
$X_1=0\klk X_{n-2}=0$ for which $f_L=0$ holds. Let $A_0$ denote
the number of different subspaces belonging to $M^{(2)}$.
Repeating this argument for all the subspaces of $M^{(2)}$ we
obtain
\begin{equation}\label{eq:estimate_C/A}
\frac{C}{A}\le
\frac{\delta^2q^{n-4}A_0}{q^{n-2}A_0}=\frac{\delta^2}{q^2}.
\end{equation}
Let us observe that $\#
M_T^{(2)}=q^n(q^n-1)(q^n-q)/\big(q^2(q^2-1)(q^2-q)\big)$.
Combining this observation with (\ref{eq:estimate_A}) we have
\begin{equation}\label{eq:estimate_D/A}
\frac{D}{A}=\frac{1}{A}(\# M_T^{(2)}-A)=\frac{1}{A}\,
\frac{q^n(q^{n-1}-1)(q^{n-1}-q)}{q^2
(q^2-1)(q^2-q)}\le\frac{4}{3q^2}.\end{equation}
Let us fix a point $x\in\fq^n$. Then there are
$E=(q^n-1)(q^n-q)/\big((q^2-1)(q^2-q)\big)$ varieties $L\in
M_T^{(2)}$ passing through $x$. This implies
\begin{equation}\label{eq:estimate_E/A}
\frac{A}{E}\le q^{n-2}.\end{equation}

Now we are ready to state our estimate on the number of rational
points of an absolutely irreducible $\fq$--hypersurface {\em
without any regularity condition}.
\begin{theorem}\label{theorem:est_pol_abs_irr}
For an absolutely irreducible $\fq$--hypersurface $H$ of $\A^n$ of
degree $\delta$ the following estimate holds:
$$ |\#(H\cap \fq^n)-q^{n-1}|\leq
(\delta-1)(\delta-2)q^{n-3/2}+5\delta^{\frac{13}{3}}q^{n-2} $$
\end{theorem}
\begin{proof}
First we observe that the theorem is obviously true if $\delta=1$.
Therefore, we may assume without loss of generality that
$\delta\ge 2$.

Let $N:=\#(H\cap\fq^n)$. With the notations introduced before, we
have:
\begin{equation}\label{sinregularidad}|N-q^{n-1}|\leq
\frac{1}{E}\bigg (\sum_{L\in M^{(2)}} |N(f_L)-q| + \sum_{L\in
M_T^{(2)}\setminus M^{(2)}} |N(f_L)-q| \bigg).\end{equation}
In order to estimate the first term of the right--hand side of
(\ref{sinregularidad}), for a plane $L \in \Pi_j$ with $j\in
\{0\klk \delta-1\}$, Lemma \ref{lemaSchmidt} implies
$|N(f_L)-q| \leq |N(f_L)-\nu(L)q|+|\nu(L)-1|q\leq \omega(q,\delta)
+ \delta^2 + jq.$
Therefore, we have:
$$\begin{array}{ccl}
\displaystyle\sum_{L\in M^{(2)}} |N(f_L)-q|& \leq &\displaystyle
\sum_{j=0}^{\delta-1}\Big(\sum_{L\in\Pi_j}\big(\omega(q,\delta)
+\delta^2+jq\big)\Big)+\sum_{L\in\Pi_{q-1}}(q^2-q)\\
%
%
& \leq &
\displaystyle\Big(\sum_{j=0}^{\delta-1}\#(\Pi_j)\Big)\big(
\omega(q,\delta) +\delta^2\big)+q\sum_{j=1}^{q-1}j\#(\Pi_j) \\
& \leq & A(\omega(q,\delta) +\delta^2) + B q + Cq(q-1).
\end{array}$$
Replacing this inequality in (\ref{sinregularidad}) and taking
into account (\ref{eq:estimate_B/A}), (\ref{eq:estimate_C/A}),
(\ref{eq:estimate_D/A}), (\ref{eq:estimate_E/A}) we obtain for
$\delta\ge 3$:
\begin{equation}\label{eq:proof_est_final_abs_irred}
\begin{array}{ccl}
|N-q^{n-1}|& \leq & \frac{1}{E} \left(A(\omega(q,\delta)
+\delta^2) + Bq + Cq(q-1)+ Dq^2\right)\\
& \leq &  \frac{A}{E}\left( \omega(q,\delta) +\delta^2 +
 \frac{B}{A}q + \frac{C}{A}q(q-1)
+\frac{D}{A}q^2\right)\\
& \leq & q^{n-2}\left(\omega(q,\delta) +\delta^2 + (2
\delta^{\frac{13}{3}} + 3 \delta^{\frac{11}{3}})\frac{4}{3}
+\delta^2+\frac{4}{3}\right)\\
& \leq & (\delta-1)(\delta-2)q ^{n-3/2} +
5\delta^{\frac{13}{3}}q^{n-2}.
\end{array}\end{equation}

For $\delta=2$, combining (\ref{sinregularidad}) with estimate
(\ref{d=2}) of the proof of Proposition \ref{prop:cota}, we obtain
\begin{equation}\label{eq:proof_est_final_abs_irred d=2}
\begin{array}{ccl}|N-q^{n-1}| &\leq &
q^{n-2}\left(\omega(q,\delta) +\delta^2 +(\frac{3}{2}\delta^4
-2\delta^3 + \frac{5}{2}\delta^2)\frac{4}{3}
+\delta^2+\frac{4}{3}\right)\\
&\leq &(\delta-1)(\delta-2)q^{n-3/2}+(2\delta^4+3\delta)q^{n-2}\\
&\leq &(\delta-1)(\delta-2)q^{n-3/2}+
5\delta^{\frac{13}{3}}q^{n-2}.
\end{array}\end{equation}
This finishes the proof of the theorem.
\end{proof}

Observe that our  estimate holds with no restriction on $q$ and
clearly improves the previous record estimate (up to the authors
knowledge), due to \cite{HuWo98}, which is only valid for $q >
cn^3\delta^5 \log ^3 \delta$:
$$|\#(H\cap \fq^n)-q^{n-1}|\leq
(\delta-1)(\delta-2)q^{n-3/2}+(2\delta^5
+\delta^2)q^{n-2}+2\delta^7q^{n-5/2}).$$

Moreover, we also improve the estimate of \cite{GhLa02a},
\cite{GhLa02}. We recall that in the case of a hypersurface the
estimate is
$$|\#(H\cap \fq^n)-q^{n-1}|\leq
(\delta-1)(\delta-2)q^{n-3/2}+12(\delta + 3)^{n+1}q^{n-2}.$$

In \cite{ChRo96}, the authors show that for $q$ sufficiently large
the following assertions hold (see \cite[Theorem 3.2 and
3.4]{ChRo96}) for any polynomial $f\in\fq[\xon]$ of degree
$\delta>0$:
\begin{enumerate}
\item[$(i)$] if $f$ is absolutely irreducible, then $\#\big(V(f)
\cap\fq^n\big) \le \delta q^{n-1}-(\delta-1)q^{n-2}$,
\item[$(ii)$] if $f$ has an absolutely irreducible nonlinear
$\fq$--definable factor, then \linebreak
$\#\big(V(f)\cap\fq^n\big) \le \delta q^{n-1}-(\delta-1)q^{n-2}$.
\end{enumerate}
Further, they ask whether the previous assertions hold for any
$q$. Although we are not able to answer this question, our
estimates provide explicit values $q_0=q_0(\delta)$ and
$q_1=q_1(\delta)$ such that $(i)$ holds for $q\ge q_0$ and $(ii)$
holds for $q\ge q_1$. Indeed, Theorem
\ref{theorem:est_pol_abs_irr} implies that we may choose
$q_0:=13\delta^{\frac{10}{3}}$ and $q_1:=9\delta^{\frac{13}{3}}$.
%
%
\subsection{An improved estimate with regularity condition.}
\label{section:lower bound}
In this section we are going to exhibit an estimate on the number
of $q$--rational points of an absolutely irreducible
$\fq$--variety which improves that of Theorem
\ref{theorem:est_pol_abs_irr} but is only valid under a certain
regularity condition.
\begin{theorem}\label{theorem:lower bound}
Let $q >15\delta^{\frac{13}{3}}$ and let $H\subset\A^n$ be an
absolutely irreducible $\fq$--hypersurface of degree $\delta$.
Then the following estimate holds:
$$ |\#(H\cap \fq^n)-q^{n-1}|\leq
(\omega(q,\delta)+5\delta^2)q^{n-2}.$$
\end{theorem}
\begin{proof}
Let $N:=\#(H\cap \fq^n)$. Since the statement of the theorem is
obviously true for $\delta=1$, we may assume without loss of
generality that $\delta\ge 2$ holds.

With notations as before, for a plane $L \in \Pi_j$ with $j>0$ it
follows by Lemma \ref{lemaSchmidt} that $|N(f_L)-q|< jq +
\omega(q,\delta)+\delta^2$ holds.
Therefore, we have 
$$\begin{array}{ccl}
|N(f_L)-Nq^{2-n}| & \geq & |N(f_L)-q|- q^{2-n}|N-q^{n-1}| \\
& \geq & jq - \omega(q,\delta)-\delta^2
-\omega(q,\delta)-5\delta^{\frac{13}{3}}\\
& \geq & \frac{1}{2}jq,
\end{array}$$
where the last inequality is valid if and only if $\frac{1}{2}jq
\geq 2q^{1/2}(\delta-1)(\delta-2) +2(\delta+1)+\delta^2
+5\delta^{\frac{13}{3}}$ holds. Hence,  our assumption on $q$
implies the validity of the inequality.

From \cite[Lemma 6]{Schmidt74} we have
$\frac{1}{4}q^2\sum_{j=1}^{q-1}j^2\#(\Pi_j)\leq \delta E q^{n-1}$,
which implies $\sum_{j=1}^{q-1}j\#(\Pi_j)\leq 4\delta Eq^{n-3}$.
Hence,
$$\begin{array}{ccl}
|N-q^{n-1}|&=&
\frac{1}{E}\Big|\sum_{L\in M^{(2)}_T}\big(N(f_L)-q\big)\Big|\leq
\frac{1}{E}\sum_{L\in M^{(2)}_T}|(N(f_L)-q)|\\
& \leq &
\frac{1}{E}\Big((A + D)
\big(\omega(q,\delta)+\delta^2\big)+\sum_{j=1}^{q-1}
2jq\#(\Pi_j)\Big)\\
& \leq & q^{n-2}\big(\omega(q,\delta)+\delta^2 +8\delta\big) \\
& \leq  & q^{n-2}\big(\omega(q,\delta)+5\delta^2\big).
\end{array}$$
This finishes the proof of the theorem.
\end{proof}

%
From this estimate we deduce the following (non--trivial) lower
bound: for $q > 1{5}\delta^{\frac{13}{3}}$, we have
$$N > q^{n-1}-(\delta-1)(\delta-2)q^{n-\frac{3}{2}}-
(5\delta^2+\delta+1)q^{n-2}.$$
Therefore, our estimate significantly improves the regularity
condition $q>10^4n^3\delta^5\vartheta^3([4\log\delta])$ of
\cite{Schmidt74}, where $[\ ]$ denotes integer part and
$\vartheta(j)$ is the $j$--th prime, and also provides a
corresponding upper bound, not given in \cite{Schmidt74}.

Let us observe that, in the setting of polynomial equation solving
over finite fields, lower bounds on the number of $q$--rational
points of a given absolutely irreducible $\fq$--hypersurface $H$,
such as those underlying Theorems \ref{theorem:est_pol_abs_irr}
and \ref{theorem:lower bound} or \cite{Schmidt74}, are typically
required in order to assure the existence of a $q$--rational point
of $H$ (see e.g. \cite{HuWo98}, \cite{HuWo99}, \cite{CaMa03}).
Indeed, from \cite{Schmidt74} one deduces that an absolutely
irreducible $\fq$--hypersurface of degree $\delta$ has a
$q$--rational point for
$q>10^4n^3\delta^5\vartheta^3([4\log\delta])$. Furthermore, from
Theorem \ref{theorem:est_pol_abs_irr} we conclude that this
condition can be improved to $q>9\delta^{\frac{13}{3}}$.
Nevertheless, the following simple argument allows us to
significantly improve the latter (compare \cite[Section
6.1]{CaMa03}):
\begin{theorem}\label{th:existenceZeros}
For $q>2\delta^4$, any absolutely irreducible $\fq$--hypersurface
of degree $\delta$ has a $q$--rational zero.
\end{theorem}
\begin{proof}
Let $H\subset\A^n$ be an absolutely irreducible
$\fq$--hypersurface of degree $\delta$, and let $f\in\fq[\xon]$ be
the defining polynomial of $H$. Since $q>2\delta^4$, from
Corollary \ref{corokaltofen} we conclude that there exists
$(\nu,\omega,\eta)\in\fq^{3n-2}$ for which
$\chi(X,Y):=f(X+\nu_1,\omega_2 X+\eta_2Y+\nu_2\klk \omega_n
X+\eta_nY+\nu_n)$ is absolutely irreducible of degree $\delta$.
Therefore, Weil's estimate (\ref{estimateWeil}) shows that $\chi$
has at least $q-(\delta-1)(\delta-2)q^{\frac{1}{2}}-\delta-1$
$q$--rational zeros. Since this quantity is a strictly positive
real number for $q>2\delta^4$, we conclude that $\chi$ has at
least one $q$--rational zero, which implies that $H$ has at least
one $q$--rational zero.
\end{proof}

Finally, we observe that in the case that the characteristic $p$
of the field $\fq$ is large enough, the estimates of Theorems
\ref{theorem:est_pol_abs_irr} and \ref{theorem:lower bound} can be
further improved, using an effective version of the first Bertini
Theorem due to S. Gao \cite{Gao03}. From \cite[Theorem 5.1]{Gao03}
we deduce the following result:
\begin{corollary}\label{coro:Gao}
Suppose that the characteristic $p$ of $\fq$ satisfies the
condition $p>2\delta^2$. Let $f \in \fq[\xon]$ be an absolutely
irreducible polynomial of degree $\delta>1$. Then there are at
most $\frac{3}{2}\delta^3\frac{q^{3n-3}}{q^3(q-1)}$ $\fq$--planes
$L\subset\A^n$ such that the restriction $f_L$ of $f$ to $L$ is
not absolutely irreducible.
\end{corollary}
With the notations of Section
\ref{subsec:hypersup_without_regularity}, from Corollary
\ref{coro:Gao} we obtain
$$\frac{B}{A}:=\frac{1}{A}\sum_{j=1}^{\delta-1}j \#\Pi_j \leq
\frac{\delta}{A}\sum_{j=1}^{\delta-1} \#\Pi_j \leq
\frac{3}{2}\frac{\delta^4}{A}\frac{q^{3n-3}}{q^3(q-1)}\leq
\frac{3}{2}\,\delta^4\frac{q^{n-2}}{q^{n-1}-1}.$$
Combining this estimate with (\ref{eq:proof_est_final_abs_irred})
of the proof of Theorem \ref{theorem:est_pol_abs_irr}, we obtain
the following estimate on the number $N$ of $q$--rational points
of any absolutely irreducible hypersurface $H\subset\A^n$ of
degree $\delta$:
$$ |N-q^{n-1}|\leq (\delta-1)(\delta-2)q^{n-3/2}+ 3\delta^4
q^{n-2}.$$
Furthermore, replacing in the proof of Theorem \ref{theorem:lower
bound} the lower bound obtained from Theorem
\ref{theorem:est_pol_abs_irr} by the one arising from the above
estimate, we obtain
$$ |N-q^{n-1}|\leq (\delta-1)(\delta-2)q^{n-3/2}+ (5\delta^2 +
\delta+1) q^{n-2}.$$
Summarizing, we have:
\begin{corollary}
\label{th:nuestroGao} Suppose that $p>2\delta^2$ holds, and let
$H\subset\A^n$ be an absolutely irreducible $\fq$--hypersurface of
degree $\delta>1$. Then the following estimate holds:
$$ |N-q^{n-1}|\leq (\delta-1)(\delta-2)q^{n-3/2}+ 3\delta^4
q^{n-2}.$$
Furthermore, if in addition we have $q> 27\delta^4$, then
$$ |N-q^{n-1}|\leq (\delta-1)(\delta-2)q^{n-3/2}+ (5\delta^2 +
\delta+1) q^{n-2}.$$
\end{corollary}
These estimates certainly improve those of Theorems
\ref{theorem:est_pol_abs_irr} and \ref{theorem:lower bound} for
$p>2\delta^2$, but fail to improve the existence result of Theorem
\ref{th:existenceZeros}. Indeed, Corollary \ref{th:nuestroGao}
does not yield a non--trivial lower bound on the number of
$q$--rational zeros of $H$ for $q\le 4\delta^4$. In fact, taking
into account that estimates like of those of Theorems
\ref{theorem:est_pol_abs_irr} and \ref{theorem:lower bound} and
Corollary \ref{th:nuestroGao} will fail to provide nontrivial
lower bounds for $q\le (\delta-1)^2(\delta-2)^2$, we conclude that
our existence result of Theorem \ref{th:existenceZeros} comes
quite close to this optimal value.
%
%
\subsection{An estimate for an arbitrary $\fq$--hypersurface.}
We finish our discussion on estimates on the number of
$q$--rational points of an $\fq$--hypersurface by considering the
case of an arbitrary $\fq$--hypersurface. Nevertheless, it must be
remarked that our estimate in the case of an $\fq$--hypersurface
without absolutely irreducible $\fq$--definable components reduces
essentially to Lemma \ref{lemma:pnai}.
\begin{theorem}
\label{theorem:estimate_hyper_gral}
Let $H\subset\A^n$ $(n\ge 2)$ be an $\fq$--hypersurface of  degree
$\delta>0$. Let $H=H_1\cup\cdots\cup H_{\sigma}\cup
H_{\sigma+1}\cup\cdots \cup H_m$ be the decomposition of $H$ into
$\fq$--irreducible components, where $H_1,\ldots,H_{\sigma}$ are
absolutely irreducible and $H_{\sigma+1},\ldots,H_m$ are not
absolutely irreducible. Let $\delta_i:=\deg H_i$ for $1\le i\le m$
and let $\Delta:=\sum_{i=1}^{\sigma}\delta_i$. Then we have the
estimate
$$|\#(H\cap \fq^n)-\sigma q^{n-1}|\leq
sign(\sigma)(\Delta-1)(\Delta-2)q^{n-\frac{3}{2}} +
(5\Delta^{\frac{13}{3}} +\delta^2/4) q^{n-2},$$
where $sign(\sigma):=0$ for $\sigma=0$ and $sign(\sigma):=1$
otherwise.
\end{theorem}
\begin{proof} 
Let $N:=\#(H \cap\fq^n)$ and $N_i:=\#(H_i \cap\fq^n)$ for $1\le
i\le m$. We have
$$ |N-\sigma q^{n-1}|\le |N-\sum_{i=1}^{\sigma}N_i| +
\sum_{i=1}^{\sigma}|N_i-q^{n-1}|.$$

For $\sigma+1 \leq i \leq m$ we have that $H_i$ is an
$\fq$--irreducible hypersurface which is not absolutely
irreducible. Therefore, Lemma \ref{lemma:pnai} implies
\begin{equation}\label{estimate:aux_1_pol}N-\displaystyle
\sum_{i=1}^{\sigma}N_i \le \displaystyle \sum_{i=\sigma+1}^m N_i <
\displaystyle q^{n-2}\sum_{i=\sigma+1}^m \delta_i^2/4 \le
q^{n-2}\delta^2/4.\end{equation}
On the other hand, we have
\begin{equation}\label{estimate:aux_2_pol} \sum_{i=1}^{\sigma} N_i-
N \le \sum_{1\leq i<j\leq \sigma}\# (H_i \cap H_j \cap \fq^n) \leq
q^{n-2}\!\!\!\!\sum_{1\leq i<j\leq \sigma}\!\!\delta_i\delta_j
\leq q^{n-2}\delta^2/4.
\end{equation}
From (\ref{estimate:aux_1_pol}) and (\ref{estimate:aux_2_pol}) we
conclude that the following estimate holds:
\begin{equation}\label{estimate:aux_3_pol}
|N-\sum_{i=1}^{\sigma}N_i|\leq q^{n-2}\delta^2/4.
\end{equation}

Since $H_i$ is an absolutely irreducible $\fq$--hypersurface of
$\A^n$ for $1\leq i \leq \sigma$, applying Theorem
\ref{theorem:est_pol_abs_irr} we obtain:
%
$$
\begin{array}{ccl}\displaystyle
\sum_{i=1}^{\sigma}|N_i-q^{n-1}|& \leq &\displaystyle
q^{n-2}\sum_{i=1}^{\sigma}\big((\delta_i-1)(\delta_i-2)q^{1/2}+ 5
\delta_i^{\frac{13}{3}}\big)\\ & \leq &
sign(\sigma)(\Delta-1)(\Delta-2)q^{n-\frac{3}{2}}+
5\Delta^{\frac{13}{3}}q^{n-2}.\end{array}
$$

Combining this estimate with (\ref{estimate:aux_3_pol}) finishes
the proof of the theorem.
%
\end{proof}
%
%
\section{From hypersurfaces to varieties.}
\label{section:reduction}
Let $V$ be an equidimensional $\fq$--variety  of dimension $r>0$
and degree $\delta$. In this section we are going to exhibit a
reduction of the problem of estimating the number of $q$--rational
points of $V$ to the hypersurface case. It is a well--known fact
that a generic linear projection morphism $\pi:\A^n\to\A^{r+1}$
induces a birational morphism which maps $V$ into a hypersurface
of $\A^{r+1}$. Our next result yields an upper bound on the degree
of the genericity condition underlying the choice of the
projection morphism $\pi$.
%
%
\begin{proposition}
\label{proposition:chow} Let $\Lambda:=(\Lambda_{ij})_{1\le i\le
r+1,1\le j\le n}$ be an $(r+1)\times n$--matrix of indeterminates,
let $\Lambda^{(i)}:=(\Lambda_{i,1}\klk \Lambda_{i,n})$ for $1\le
i\le r+1$, and let $\Gamma:=(\Gamma_1\klk \Gamma_{r+1})$ be an
$(r+1)$--dimensional vector of indeterminates. Let $X:=(\xon)$ and
let $\widetilde{Y}:=\Lambda X+\Gamma$. Then there exists a nonzero
polynomial $G\in\cfq[\Lambda,\Gamma]$  of degree at most $
2(r+1)\delta^2$ such that for any
$(\lambda,\gamma)\in\A^{(r+1)n}\times\A^{r+1}$ with
$G(\lambda,\gamma)\not=0$ the following conditions are satisfied:
\begin{enumerate}
  \item[$(i)$] Let $Y:=\lambda X+\gamma:=(\yo{r+1})$. Then the projection
  morphism \linebreak $\pi:V\to\A^{r}$ defined by $\yo{r}$ is a finite morphism.
  \item[$(ii)$] The linear form $Y_{r+1}$ induces a primitive element
  of the integral ring extension $\cfq[\yo{r}]\hookrightarrow\cfq[V]$,
  i.e. the degree of its minimal integral dependence equation in
  $\cfq[\yo{r}]$ equals the rank of $\cfq[V]$ as (free)
  $\cfq[\yo{r}]$--module.
\end{enumerate}
\end{proposition}

\begin{proof}
Let us consider the following morphism of algebraic varieties:
\[\begin{array}{crcl}
  \Phi:& \A^{(r+1)n}\times\A^{r+1}\times V &
  \to & \A^{(r+1)n}\times\A^{r+1}\times \A^{r+1} \\
       & (\lambda,\gamma,x) & \mapsto & (\lambda,\gamma,\lambda
       x+\gamma).
\end{array}\]
Using standard facts about Chow forms (see e.g.
\cite{Shafarevich84}, \cite{KrPaSo01}), we conclude that
$\overline{Im(\Phi)}$ is a hypersurface of
$\A^{(r+1)n}\times\A^{r+1}\times \A^{r+1}$, defined by a
polynomial $P\in\cfq[\Lambda,\Gamma,\widetilde{Y}_1\klk
\widetilde{Y}_{r+1} ]$ which satisfies the following estimates:
\begin{itemize}
  \item $\deg_{\widetilde{Y}_1\klk
\widetilde{Y}_{r+1}}P= \deg_{\widetilde{Y}_{r+1}}P= \delta$,
  \item $\deg_{\Lambda^{(i)}\!,\,\Gamma_i}P\le \delta$
for $1\le i\le r+1$.
\end{itemize}

Let $G_1\in\cfq[\Lambda,\Gamma]$ be the (nonzero) coefficient of
the monomial $\widetilde{Y}_{r+1}^\delta$ in the polynomial $P$,
considering $P$ as an element of
$\cfq[\Lambda,\Gamma][\widetilde{Y}_1\klk \widetilde{Y}_{r+1}]$.
We have $\deg G_1\le (r+1)\delta$. Let
$\widetilde{G}_1\in\cfq[\Lambda^{\!(i)},\Gamma_i:1\le i\le r]$ be
the coefficient of a nonzero monomial of the polynomial $G_1$,
considering $G_1$ as an element of
$\cfq[\Lambda^{\!(i)},\Gamma_i:1\le i\le
r][\Lambda^{(r+1)},\Gamma_{r+1}]$.

Let $(\lambda^*,\gamma^*)\in\A^{rn}\times\A^{r}$ be any point
satisfying the condition
$\widetilde{G}_1(\lambda^*,\gamma^*)\not=0$ and let $Y:=(Y_1\klk
Y_r):=\lambda^* X+\gamma^*$. We claim that condition ($i$) of the
statement of Proposition \ref{proposition:chow} holds. Indeed,
since $G_1^*:=G_1(\lambda^*,\gamma^*,\Lambda^{(r+1)},
\Gamma_{r+1})$ is a nonzero element of
$\cfq[\Lambda^{(r+1)},\Gamma_{r+1}]$, we deduce that there exist
$n$ $\cfq$--linearly independent vectors $w_1\klk w_n\in\cfq^n$
and $a_1\klk a_n\in\cfq$ such that $G_1^*(w_k,a_k)\not=0$ holds
for $1\le k\le n$. Let $\ell_k:=w_k X+a_k$ for $1\le k\le n$. By
construction, for $1\le k\le n$ the polynomial
$P(\lambda^*,\gamma^*,w_k,a_k,\yo{r},\ell_k)$ is an integral
dependence equation over $\cfq[\yo{r}]$ for the coordinate
function induced by $\ell_k$ in $\cfq[V]$. Since
$\cfq[\ell_1\klk\ell_n]= \cfq[\xon]$ we conclude that condition
($i$) holds.

Furthermore, since $\cfq[\Lambda,\Gamma,\widetilde{Y}_1\klk
\widetilde{Y}_{r+1}]/ (P)$ is a reduced $\cfq$--algebra and $\cfq$
is a perfect field, using \cite[Proposition 27.G]{Matsumura80} we
see that the (zero--dimensional)
$\cfq(\Lambda,\Gamma,\widetilde{Y}_1\klk\widetilde{Y}_{r})$--algebra
$\cfq(\Lambda,\Gamma,\widetilde{Y}_1\klk\widetilde{Y}_{r})
[\widetilde{Y}_{r+1}]/ (P)$ is reduced. This implies that $P$ is a
separable element of
$\cfq(\Lambda,\Gamma,\widetilde{Y}_1\klk\widetilde{Y}_{r})
[\widetilde{Y}_{r+1}]$, and hence $P$ and $\partial P/\partial
\widetilde{Y}_{r+1}$ are coprime in
$\cfq(\Lambda,\Gamma,\widetilde{Y}_1\klk\widetilde{Y}_r)
[\widetilde{Y}_{r+1}]$. Then the discriminant
$\rho:=Res_{\widetilde{Y}_{r+1}}(P,\partial P/\partial
\widetilde{Y}_{r+1})$ of $P$ with respect to $\widetilde{Y}_{r+1}$
is a nonzero element of
$\cfq[\Lambda,\Gamma,\widetilde{Y}_1\klk\widetilde{Y}_r]$ which
satisfies the following degree estimates:
\begin{itemize}
  \item $\deg_{\widetilde{Y}_1\klk\widetilde{Y}_r}\rho\le
(2\delta-1)\delta$,
  \item $\deg_{\Lambda^{(i)}\!,\,\Gamma_i}\rho\le
(2\delta-1)\delta$ for $1\le i\le r+1$.
\end{itemize}

Let $\rho_1\in\cfq[\Lambda,\Gamma]$ be a (nonzero) coefficient of
a monomial of $\rho$, considering $\rho$ as an element of
$\cfq[\Lambda,\Gamma][\widetilde{Y}_1\klk \widetilde{Y}_r]$ and
let $G:=\rho_1\,\widetilde{G}_1$. Observe that $\deg G\le
2(r+1)\delta^2$ holds.

Let $(\lambda,\gamma)\in\A^{(r+1)n}\times \A^{r+1}$ satisfy the
condition $G(\lambda,\gamma)\not=0$, let
$(\lambda^*,\gamma^*)\in\A^{rn}\times\A^r$ be the first $r$ rows
of $(\lambda,\gamma)$ and let $Y:=(\yo{r+1})=\lambda X+\gamma$. It
is clear that condition ($i$) holds.

We are going to prove that condition ($ii$) holds. For this
purpose, let $\rho^*$ be the polynomial obtained from $\rho$ by
specializing $\Lambda^{\!(i)},\Gamma_i$ ($1\le i\le r$) into the
value $(\lambda^*,\gamma^*)$. Then $\rho^*$ is a nonzero
polynomial of $\cfq[\Lambda^{(r+1)},\Gamma_{r+1},\yo{r}]$ which
equals the discriminant of
$P(\lambda^*,\Lambda^{(r+1)},\gamma^*,\Gamma_{r+1},\yo{r},
\widetilde{Y}_{r+1})$ with respect to $\widetilde{Y}_{r+1}$.

Let $\xi_1\klk \xi_n$ be the coordinate functions of $V$ induced
by $\xon$, let $\zeta_i:=\sum_{j=1}^n\lambda_{i,j}\xi_j$ for $1\le
i\le r$ and let $\widehat{Y}_{r+1}:=\sum_{j=1}^n
\Lambda_{r+1,j}\xi_j$. From the properties of the Chow form of $V$
we conclude that the identity
\begin{equation}
\label{equation:chow1}
\begin{array}{ccl}
0 &=&P(\lambda^*\!,\Lambda^{\!(r+1)}\!,\gamma^*\!,
\Gamma_{r+1},\zeta_1\klk \zeta_r,\widehat{Y}_{r+1})\\
&=&P(\lambda^*\!,\Lambda^{\!(r+1)}\!,\gamma^*\!,
\Gamma_{r+1},\zeta_1,\ldots,\zeta_r,\Lambda_{r+1,1}\xi_1\!\plp
\!\Lambda_{r+1,n}\xi_n)
\end{array}
\end{equation}
holds in $\cfq[\Lambda^{(r+1)},\Gamma_{r+1}]\otimes_{\cfq}
\cfq[V]$.
Following e.g. \cite{AlBeRoWo96} or \cite{Rouillier97}, from
(\ref{equation:chow1}) one deduces that for $1\le k\le n$ the
identity
\begin{equation}
\label{equation:chow2}
\begin{array}{c}
(\partial P/\partial \widetilde{Y}_{r+1})
(\lambda^*,\Lambda^{(r+1)},\gamma^*,\Gamma_{r+1},\zeta_1 \klk
\zeta_r,\widehat{Y}_{r+1})\xi_k+ \qquad\qquad\quad\\
\qquad\qquad\qquad+(\partial P/\partial \Lambda_{k})
(\lambda^*,\Lambda^{(r+1)},\gamma^*,\Gamma_{r+1},\zeta_1\klk
\zeta_r,\widehat{Y}_{r+1})=0
\end{array}
\end{equation}
holds in $\cfq[\Lambda^{(r+1)},\Gamma_{r+1}]\otimes_{\cfq}
\cfq[V]$. Since $\rho^*(\Lambda^{\!(r+1)},\Gamma_{r+1},Y_1\klk
Y_r)$ is the discriminant of
$P(\lambda^*,\Lambda^{(r+1)},\gamma^*,\Gamma_{r+1},\yo{r},
\widetilde{Y}_{r+1})$ with respect to $\widetilde{Y}_{r+1}$, it
can be expressed as a linear combination of
$P(\lambda^*,\Lambda^{(r+1)},\gamma^*,\Gamma_{r+1},\yo{r},
\widetilde{Y}_{r+1})$ and $(\partial P/\partial
\widetilde{Y}_{r+1}) (\lambda^*, \Lambda^{(r+1)},
\gamma^*,\Gamma_{r+1},\yo{r},\widetilde{Y}_{r+1})$. Combining this
observation with (\ref{equation:chow1}) and (\ref{equation:chow2})
we conclude that
\begin{equation}
\label{equation:chow3}
\begin{array}{c}
\rho^*(\Lambda^{\!(r+1)},\Gamma_{r+1},\zeta_1\klk \zeta_r)\xi_k+
P_k (\Lambda^{\!(r+1)},\Gamma_{r+1},\zeta_1\klk
\zeta_r,\widehat{Y}_{r+1})=0
\end{array}\end{equation}
holds, where $P_k$ is a nonzero element of
$\cfq[\Lambda^{\!(r+1)},\Gamma_{r+1},\om Z {r+1}]$ for $1\le i\le
n$. Specializing identity (\ref{equation:chow3}) into the values
$\Lambda_{r+1,j}:=\lambda_{r+1,j}$ ($1\le j\le n$) and
$\Gamma_{r+1}=\gamma_{r+1}$ for $1\le k\le n$ we conclude that
$Y_{r+1}$ induces a primitive element of the $\cfq$--algebra
extension $\cfq(Y_1\klk Y_r)\hookrightarrow\cfq(Y_1\klk
Y_r)\otimes_{\cfq}\cfq[V]$.

Condition ($i$) implies that $\cfq[V]$ is a finite free
$\cfq[Y_1\klk Y_r]$--module and hence $\cfq(Y_1\klk
Y_r)\otimes_{\cfq}\cfq[V]$ is a finite--dimensional $\cfq(Y_1\klk
Y_r)$--vector space. Furthermore, the dimension of $\cfq(Y_1\klk
Y_r)\otimes_{\cfq}\cfq[V]$ as $\cfq(Y_1\klk Y_r)$--vector space
equals the rank of $\cfq[V]$ as $\cfq[Y_1\klk Y_r]$--module. On
the other hand, since $\cfq[Y_1\klk Y_r]$ is integrally closed we
have that the minimal dependence equation of any element of
$f\in\cfq[V]$ over $\cfq(Y_1\klk Y_r)$ equals the minimal integral
dependence of $f$ over $\cfq[Y_1\klk Y_r]$ (see e.g. \cite[Lemma
II.2.15]{Kunz85}). Combining this remark with the fact that
$Y_{r+1}$ induces a primitive element of the $\cfq$--algebra
extension $\cfq(Y_1\klk Y_r)\hookrightarrow\cfq(Y_1\klk
Y_r)\otimes_{\cfq}\cfq[V]$ we conclude that $Y_{r+1}$ also induces
a primitive element of the $\cfq$--algebra extension $\cfq[Y_1\klk
Y_r]\hookrightarrow\cfq[V]$. This shows condition ($ii$) and
finishes the proof of the proposition.
\end{proof}
%
%

From Proposition \ref{proposition:chow} we easily deduce that $V$
is birationally equivalent to an $\cfq$--hypersurface
$H\subset\A^{r+1}$ of degree $\delta$, namely the image of $V$
under the projection defined by linear forms $Y:=(Y_1\klk
Y_{r+1})=\lambda X+\gamma$ with $G(\lambda,\gamma)\not=0$, where
$G$ is the polynomial of the statement of Proposition
\ref{proposition:chow} (compare Proposition
\ref{proposition:morph_bir} below).

We would like to estimate the number of $q$--rational points of
the variety $V$ in terms of that of the hypersurface $H$, but
``good" estimates on the number of $q$--rational points of $H$ are
not available if $H$ is not an $\fq$--variety. Let us observe that
$H$ is an $\fq$--variety if the linear forms $Y_1\klk Y_{r+1}$
belong to $\fq[\xon]$ (see e.g. \cite{Kunz85}). In order to ensure
that there exist linear forms $Y_1\klk Y_{r+1}\in\fq[\xon]$
satisfying conditions ($i$) and ($ii$) of Proposition
\ref{proposition:chow} we have the following result:

\begin{corollary} \label{corollary:Zippel-Schwartz} Let notations
and assumptions be as in Proposition \ref{proposition:chow}. If $q
> 2(r+1)\delta^2$, there exists an element $(\lambda,\gamma)\in
\fq^{(r+1)\times n}\times\fq^{r+1}$ satisfying the condition
$G(\lambda,\gamma)\not=0$, where $G$ is the polynomial of the
statement of Proposition \ref{proposition:chow}.
\end{corollary}
\begin{proof} Let $V(G):=\{(\lambda,\gamma)\in\A^{(r+1)n}\times
\A^{r+1}:G(\lambda,\gamma)=0\}$. Taking into account the upper
bound of Lemma \ref{lemma:hesch2}
$$\#\big(V(G)\cap(\fq^{(r+1) n}\times \fq^{r+1})\big)\le
2(r+1)\delta^2q^{(r+1)(n+1)-1},$$
we immediately deduce the statement of Corollary
\ref{corollary:Zippel-Schwartz}.
\end{proof}

From now on, we shall assume that the condition $q>
2(r+1)\delta^2$ holds. Let $(\lambda,\gamma)\in\fq^{(r+1)n}\times
\fq^{r+1}$ satisfy $G(\lambda,\gamma)\not=0$, let $Y=(Y_1\klk
Y_{r+1}):=\lambda Y+\gamma$ and let us consider the following
$\fq$--definable morphism of $\fq$--varieties:
$$\begin{array}{ccrcl}
  \pi & : & V & \rightarrow & \A^{r+1} \\
   &  & x & \mapsto &
   \big(Y_1(x),\ldots,Y_{r+1}(x)\big).
\end{array}$$
Then the set $W:=\pi(V)$ is an $\fq$--hypersurface. This
hypersurface is defined by a polynomial $h\in \fq[\yo{r+1}]$,
which is a separable monic element of the polynomial ring
$\fq[\yo{r}][Y_{r+1}]$.

Let $V_1\subset\A^n$ and $W_1\subset\A^{r+1}$ be the following
$\fq$-varieties: $$\begin{array}{rl}
  V_1 := & \{x\in \A^n:
  (\partial h/\partial Y_{r+1})(Y_1(x),\ldots,Y_{r+1}(x))=0\}, \\
  W_1 := & \{y\in \A^{r+1}:
  (\partial h/\partial Y_{r+1})
 (y)=0\}.
\end{array}$$

Our following result shows that the variety $V$ is birationally
equivalent to the hypersurface $W\subset\A^{r+1}$.
%
\begin{proposition}\label{proposition:morph_bir} Let $q>2(r+1)\delta^2$.
Then ${\pi}|_{V\setminus V_1}: V \setminus V_1 \rightarrow W
\setminus W_1$ is an isormorphism of\ $\fq$--Zariski open sets.
\end{proposition}
\begin{proof}
Let us observe that $\pi(V\setminus V_1)\subset W \setminus W_1$.
Then ${\pi}|_{V\setminus V_1}: V \setminus V_1 \rightarrow W
\setminus W_1$ is a well--defined $\fq$--definable morphism.

We claim that $\pi$ is an injective mapping. Indeed, specializing
identity (\ref{equation:chow2}) of the proof of Proposition
\ref{proposition:chow} into the values
$\Lambda_{r+1,j}:=\lambda_{r+1,j}$ ($1\le j\le n$) and
$\Gamma_{r+1}=\gamma_{r+1}$ we deduce that there exist polynomials
$v_1\klk v_n\in\fq[Y_1\klk Y_{r+1}]$ such that for $1\le i\le n$
the following identity holds:
\begin{equation}\label{def:v_i}v_i(Y_1\klk Y_{r+1})-X_i
\cdot(\partial h/\partial Y_{r+1})(Y_1\klk Y_{r+1})\equiv 0 \mbox{
mod $I(V)$ }.
\end{equation}
Let $x:=(x_1\klk x_n),x':=(x_1'\klk x_n')\in V \setminus V_1$
satisfy $\pi(x)=\pi(x')$. We have $Y_i(x)=Y_i(x')$ for $1\le i\le
r+1$. Then from (\ref{def:v_i}) we conclude that $x_i=x_i'$ for
$1\le i\le n$, which shows our claim.\smallskip

Now we show that $\pi|_{V\setminus V_1}:V\setminus V_1\to
W\setminus W_1$ is a surjective mapping. Let $h_0:=\partial
h/\partial Y_{r+1}$. Let be given an arbitrary element
$y:=(y_1\klk y_{r+1})$ of $W\setminus W_1$, and let
$$x:=\big((v_1/h_0)(y),\ldots,(v_{n}/h_0)(y)\big).$$ We claim that
$x$ belongs to $V\setminus V_1$. Indeed, let $f$ be an arbitrary
element of the ideal $I(V)$ and let
$\tilde{f}:=\left(h_0(Y_1,\ldots,Y_{r+1})\right)^N f$, where
$N:=\deg f$. Then there exists $g \in \F\!_q[Z_1,\ldots,Z_{n+1}]$
such that $\tilde{f}= g(h_0X_1,\ldots,h_0X_n,h_0)$ holds. Since
$\tilde{f} \in I(V)$,  for any $z \in V$ we have $\tilde{f}(z)=0$
and hence from identity (\ref{def:v_i}) we conclude that
$g(v_1,\ldots,v_n,h_0)\big(Y_1(z)\klk Y_{r+1}(z)\big)=0$ holds.
This shows that $h$ divides $\hat{f}:=g(v_1,\ldots,v_n,h_0)$ in
$\fq[\yo{r+1}]$ and therefore $\hat{f}(y)=h_0(y)^N f(x)=0$ holds.
Taking into account that $h_0(y)\neq 0$ we conclude that $f(x)=0$
holds, i.e. $x \in V\setminus V_1$.

In order to finish the proof of the surjectivity of $\pi$ there
remains to prove that $\pi(x)=y$ holds. For this purpose, we
observe that identity (\ref{def:v_i}) shows that any $z \in V$
satisfies
$$Y_{i}(z)h_0\big(Y_1(z)\klk Y_{r+1}(z)\big)-\sum_{j=1}^n
\lambda_{i,\,j}\,v_j\big(Y_1(z)\klk Y_{r+1}(z)\big)=0$$
for $1\le i\le r+1$. Then $h$ divides the polynomial
$Y_ih_0-\sum_{j=1}^n\!\lambda_{i,j}v_j$ in $\F\!_q[Y_1,\ldots
Y_{r+1}]$, which implies $y_i=\sum_{j=1}^n
\lambda_{i,\,j}(v_j/h_0)(y)=\sum_{j=1}^n \lambda_{i,\,j}\,x_j$ for
$1\le i\le r+1$. This proves that $\pi(x)=y$ holds.\smallskip

Finally we show that $\pi|_{V\setminus V_1}:V\setminus V_1\to
W\setminus W_1$ is an isomorphism. Let $$\begin{array}{ccrcl}
  \phi & : & W\setminus W_1 & \rightarrow & V\setminus V_1 \\
   &  & y & \mapsto & \big((v_1/h_0)(y),\ldots,(v_{n}/h_0)(y)\big).
\end{array}$$
Our previous discussion shows that $\phi$ is a well--defined
$\fq$--definable morphism. Furthermore, our arguments above show
that $\pi\circ\phi$ is the identity mapping of $W\setminus W_1$.
This finishes the proof of the proposition.
\end{proof}

From Proposition \ref{proposition:morph_bir} we immediately
conclude that the $\fq$--Zariski open sets $V\setminus V_1$ and
$W\setminus W_1$ have the same number of $q$--rational points.
%
%
\section{Estimates for an $\fq$--variety.}\label{section:variety}
In this section we exhibit explicit estimates on the number of
$q$--rational points of an $\fq$--variety. For this purpose, we
are going to apply the reduction to the hypersurface case of
Section \ref{section:reduction}, together with the estimates for
hypersurfaces of Section \ref{section:hypersurface}. We start with
the case of an absolutely irreducible $\fq$--variety.
\begin{theorem}\label{theorem:est_final_abs_irr}
Let $V\subset \A^n$ be an absolutely irreducible $\fq$--variety of
dimension $r>0$ and degree $\delta$. If $q>2(r+1){\delta}^2$, then
the following estimate holds:
\begin{equation}\label{equation:est_final_abs_irr}|\#(V\cap {
\F\!_q}\!^n) -q^r|<
(\delta-1)(\delta-2)q^{r-\frac{1}{2}}+5\delta^{\frac{13}{3}}q^{r-1}.
\end{equation}
\end{theorem}
\begin{proof}
First we observe that the theorem is obviously true in the cases
$n=1$ or $\delta=1$, and follows from Weil's estimate
(\ref{estimateWeil}) in the case $n=2$. Therefore, we may assume
without loss of generality that $n\ge 3$ and $\delta\ge 2$ hold.

Since the condition $q>2(r+1){\delta}^2$ holds, from Corollary
\ref{corollary:Zippel-Schwartz} we deduce that there exist linear
forms $Y_1\klk Y_{r+1}\in\fq[\xon]$ satisfying conditions ($i$)
and ($ii$) of the statement of Proposition \ref{proposition:chow}.
Therefore, from Proposition \ref{proposition:morph_bir} we have
$$|\#(V\cap\fq^n)-q^r|\le |\#(W\cap\fq^{r+1})-q^r| + \#(V\cap
V_1\cap\fq^n) + \#(W\cap W_1\cap \fq^{r+1} ),$$
where $V_1\subset\A^n$, $W\subset\A^{r+1}$ and
$W_1\subset\A^{r+1}$ are the $\fq$--hypersurfaces defined by the
polynomials $(\partial h/\partial
Y_{r+1})\big(Y_1(X),\dots\!,\!Y_{r+1}(X) \big)\in\fq[\xon]$,
$h\in\fq[Y_1,\dots\!,\!Y_{r+1}]$ and $(\partial h/\partial
Y_{r+1})\in\fq[Y_1\klk Y_{r+1}]$ respectively.

From the B\'ezout inequality (\ref{equation:Bezout}) and Lemma
\ref{lemma:hesch2} we deduce the upper bounds:
\begin{equation}\label{equation:aux_est_1}\begin{array}{rcl} \#
(V\cap V_1\cap\fq^n)&\le & \delta(\delta-1)q^{r-1}, \\ \#(W\cap
W_1\cap \fq^{r+1} )& \le & \delta(\delta-1)q^{r-1}.
\end{array}\end{equation}
On the other hand, we observe that $W$ is an absolutely
irreducible $\fq$--variety of dimension $r>0$ and degree
$\delta>0$. Therefore, applying estimate
(\ref{eq:proof_est_final_abs_irred}) of the proof of Theorem
\ref{theorem:est_pol_abs_irr} we obtain
$$\begin{array}{c} |\#(W\cap \fq^{r+1})-q^{r}|\le
(\delta-1)(\delta-2)q^{r-\frac{1}{2}}+ \Big(\frac{8}{3}
\delta^{\frac{13}{3}} + 4 \delta^{\frac{11}{3}} +2\delta^2+\delta+
\frac{7}{3}\Big)q^{r-1}.\end{array}$$
This estimate, together with (\ref{equation:aux_est_1}),
immediately implies the statement of the theorem for $\delta \geq
3$. For $\delta=2$ we combine the above estimate with
(\ref{equation:aux_est_1}) and (\ref{eq:proof_est_final_abs_irred
d=2}), which yields the estimate of the statement of the theorem.
This finishes the proof.
\end{proof}

Furthermore, if we estimate $|\#(W\cap \fq^{r+1})-q^{r}|$ using
Theorem \ref{theorem:lower bound} instead of Theorem
\ref{theorem:est_pol_abs_irr}, we obtain the following result:
\begin{corollary}
\label{coro:varietyWithRegularity} Let $V\subset \A^n$ be an
absolutely irreducible $\fq$--variety of dimension $r>0$ and
degree $\delta$. If
$q>\max\{2(r+1){\delta}^2,1{5}\delta^{\frac{13}{3}}\}$, then the
following estimate holds:
$$|\#(V\cap \fq^{r+1})-q^{r}|\le
(\delta-1)(\delta-2)q^{r-\frac{1}{2}}+7\delta^2q^{r-1}.$$
\end{corollary}
Finally, if the characteristic $p$ of $\fq$ is greater than
$2\delta^2$, from Corollary \ref{th:nuestroGao} we obtain:
\begin{corollary} Let $V\subset \A^n$ be an absolutely irreducible
$\fq$--variety of dimension $r>0$ and degree $\delta$. If
$p>2\delta^2$ and $q>2(r+1){\delta}^2$ we have: $$|\#(V\cap
\fq^{r+1})-q^{r}|\le
(\delta-1)(\delta-2)q^{r-\frac{1}{2}}+4\delta^4q^{r-1}.$$ If in
addition  $q>25\delta^4$, then the following estimate holds:
$$|\#(V\cap \fq^{r+1})-q^{r}|\le
(\delta-1)(\delta-2)q^{r-\frac{1}{2}}+7\delta^2q^{r-1}.$$
\end{corollary}

The estimate of Theorem \ref{theorem:est_final_abs_irr} yields a
nontrivial lower bound on the number of $q$--rational points of an
absolutely irreducible $\fq$--variety $V$ of dimension $r>0$ and
degree $\delta$, implying thus the existence of a $q$--rational
point of $V$, for
$q>\max\{2(r+1)\delta^2,9\delta^{\frac{13}{3}}\}$. Nevertheless,
similarly to Theorem \ref{th:existenceZeros}, the following simple
argument allows us to obtain the following improved existence
result:
\begin{corollary} For $q>\max\{2(r+1)\delta^2,2\delta^4\}$, any
absolutely irreducible $\fq$--variety $V$ of dimension $r>0$ and
degree $\delta$ has a $q$--rational point.
\end{corollary}
\begin{proof}
Since $q>2(r+1)\delta^2$ holds, from Corollary
\ref{corollary:Zippel-Schwartz} we conclude that there exist
linear forms $Y_1\klk Y_{r+1}\in\fq[\xon]$ satisfying the
conditions of Proposition \ref{proposition:chow}. Let
$h\in\fq[\yo{r+1}]$ denote the defining polynomial of the
absolutely irreducible $\fq$--hypersurface $W\subset\A^{r+1}$
defined by the image of the linear projection of $V$ induced by
$Y_1\klk Y_{r+1}$. From the condition $q>2\delta^4$ we conclude
that there exists an $\fq$--plane $L\subset\A^{r+1}$ for which
$W\cap L$ is an absolutely irreducible $\fq$--curve of $\A^{r+1}$.
Hence, Weil's estimate (\ref{estimateWeil}) shows that $\#(W\cap
L\cap\fq^{r+1})>q-(\delta-1)(\delta-2)q^{\frac{1}{2}}-\delta-1$
holds. Furthermore, from the B\'ezout inequality we deduce that
$\#(W\cap L\cap V(\partial h/\partial Y_{r+1}))\le
\delta(\delta-1)$ holds, which implies $\#\Big(\big(W\setminus
V(\partial h/\partial Y_{r+1})\big)\cap L\cap\fq^{r+1}\Big)>
q-(\delta-1)(\delta-2)q^{\frac{1}{2}}-\delta^2-1$. Since this
quantity is strictly positive for $q>2\delta^4$, it follows that
there exists a $q$--rational point of $W\setminus V(\partial
h/\partial Y_{r+1})$. Combining this with Proposition
\ref{proposition:morph_bir} we conclude that there exists a
$q$--rational point of $V$.
\end{proof}
\subsection{An estimate for an arbitrary $\fq$--variety.}
Now we are going to estimate the number of $q$--rational points of
an arbitrary $\fq$--variety $V$ of dimension $r>0$ and degree
$\delta$. Let $V= V_1\cup \cdots \cup V_{m}$ be the decomposition
of $V$ into $\fq$--irreducible components, and suppose that the
numbering is such that $V_i$ is absolutely irreducible of
dimension $r>0$ for $1\le i\le \sigma$, absolutely irreducible of
dimension at most $r-1$ for $\sigma+1\le i\le \rho$ and not
absolutely irreducible for $\rho+1\le i\le m$.

For $1\le i\le m$, let $N_i:=\#(V_i \cap \fq^n)$ and denote by
$\delta_i$ the degree of $V_i$. Finally, let
$\Delta:=\sum_{i=1}^{\sigma}\delta_i$ and $N:=\#(V \cap \fq^n)$.
We have the following result:
\begin{theorem}\label{theorem:est_final_variety}
With notations and assumptions as above, if $q
>2(r+1)\delta^2$ the number $N$ of $q$--rational points of the
variety $V$ satisfies the following estimate:
\begin{equation}\label{estimate:final}|N-\sigma q^r|\le
sign(\sigma)(\Delta-1)(\Delta-2)q^{r-1/2}+
\big(5\Delta^{\frac{13}{3}}+\delta^2\big)q^{r-1},
\end{equation}
where $sign(\sigma):=0$ for $\sigma=0$ and $sign(\sigma):=1$
otherwise.
\end{theorem}
\begin{proof}
We have $|N-\sigma q^r| \leq \sum_{i=1}^{\sigma}|N_i-q^r|+
|N-\sum_{i=1}^{\sigma}N_i|$.

From Theorem \ref{theorem:est_final_abs_irr} we obtain:
\begin{equation}\label{aux:second}
\begin{array}{ccl}
\displaystyle \sum_{i=1}^{\sigma}|N_i-q^{r}|& \leq & \displaystyle
\sum_{i=1}^{\sigma}\big( (\delta_i-1)(\delta_i -2)q^{r-1/2}+
5\delta_i^{\frac{13}{3}}q^{r-1}\big)\\ & \le &
sign(\sigma)(\Delta-1)(\Delta-2)q^{r-1/2} +
5\Delta^{\frac{13}{3}}q^{r-1}.
\end{array}
\end{equation}

Now we estimate the term $ |N-\sum_{i=1}^{ \sigma}N_i|$. Let
$\sigma+1\le i\le \rho$. Then $V_i$ is an $\fq$--variety of
dimension at most $r-1$ and degree $\delta_i$, and Lemma
\ref{lemma:hesch2} implies $N_i \leq \delta_iq^{r-1}$. On the
other hand, for $\rho+1\le i\le m$ we have that $V_i$ is
$\fq$--irreducible and not absolutely irreducible, and Lemma
\ref{lemma:pnai} shows that $N_i \leq \delta_i^2q^{r-1}/4$ holds.
Then we have
\begin{equation}\label{estimate:11}N-\displaystyle
\sum_{i=1}^{\sigma}N_i  \leq  \displaystyle \sum_{i=\sigma+1}^m
N_i \le   q^{r-1}\!\!\!\sum_{i=\sigma+1}^m \delta_i^2 \le
\delta^2q^{r-1}.
\end{equation}

On the other hand, Lemma \ref{lemma:hesch2} implies
\begin{equation}\label{estimate:22}
\sum_{i=1}^{\sigma}N_i-N  \le \sum_{1\leq i<j\leq
\sigma}\!\!\!\!\# (V_i \cap V_j \cap \fq^n) \le q^{r-1}
\!\!\!\!\!\!\sum_{1\leq i<j\leq \sigma}\!\!\!\!\delta_i \delta_j
\le \delta^2q^{r-1}/2.
\end{equation}
From estimates (\ref{estimate:11}) and (\ref{estimate:22}) we
conclude that $|N-\sum_{i=1}^{\sigma}N_i|\leq \delta^2q^{r-1} $
holds. Combining this estimate with (\ref{aux:second}) finishes
the proof of the theorem.
\end{proof}

{\footnotesize

}

\end{document}